\documentclass{article}
\usepackage{bbm}
\usepackage[title,titletoc]{appendix}
\usepackage{graphicx}
\usepackage{amsmath,amssymb,amsthm,amsfonts}
\usepackage{color}
\usepackage{amsfonts}
\usepackage{empheq}
\usepackage{indentfirst}
\usepackage{cite}
\usepackage{mathrsfs}
\usepackage[colorlinks=true, linkcolor=blue, urlcolor=blue]{hyperref}
\allowdisplaybreaks[4]
\usepackage[pagewise]{lineno}
\newtheorem{thm}{Theorem}[section]
\newtheorem{lem}{Lemma}[section]

\theoremstyle{definition}

\newtheorem{rem}{Remark}[section]

\numberwithin{equation}{section}

\DeclareMathSymbol{\C}{\mathalpha}{AMSb}{"43} \topmargin-.1in
\def\s{\hat{\tilde{S}}}
\def\f{\frac}
\def\F{\displaystyle\frac}
\def\dis{\displaystyle}
\def\I{\dis\int}
\def\i{\int}
\def\Io{\dis\int_0^L}
\def\io{\int_0^L}

\newcommand{\beq}{\begin{equation}}
\newcommand{\eeq}{\end{equation}}
\newcommand{\bea}{\begin{array}}
\newcommand{\eea}{\end{array}}

\newcommand{\eps}{\varepsilon}
\newcommand{\kap}{\kappa}
\newcommand{\lam}{\lambda}
\newcommand{\Lam}{\Lambda}

\newcommand{\sig}{\sigma}
\newcommand{\rr}{\mathcal{R}}
\newcommand{\gam}{\gamma}

\newcommand{\del}{\delta}
\newcommand{\al}{\alpha}
\newcommand{\be}{\beta}
\def\th{\theta}
\newcommand{\pa}{\partial}
\def\dx{\mathrm{d}x}
\def\dy{\mathrm{d}y}

\def\ni{\noindent}
\def\proof{{\ni\bf Proof:\quad}}
\def\proofend{{\hfill$\Box$}\\}
\def\one{{\ni\bf Proof of Theorem \ref{properties of R0}: }}
\def\oneplus{{\ni\bf Proof of Theorem \ref{lim of lam1}: }}
\def\two{{\ni\bf Proof of Theorem \ref{global stablity of ee}: }}
\def\three{{\ni\bf Proof of Theorem \ref{ee q-infty}: }}
\def\four{{\ni\bf Proof of Theorem \ref{ee ds-0}: }}
\def\five{{\ni\bf Proof of Theorem \ref{ee di-0}: }}
\def\six{{\ni\bf Proof of Theorem \ref{ee ds-infty}: }}
\def\seven{{\ni\bf Proof of Theorem \ref{ee di-infty}: }}
\def\eight{{\ni\bf Proof of Theorem \ref{ee m-infty}: }}

\title{Qualitative analysis on a spatial SIS  epidemic model with linear source in advective environments II saturated incidence}

\author{Qi Wang
\thanks{College of Science, University of Shanghai for Science and Technology,
Shanghai 200093, P.R. China.
Email: \texttt{qwang@usst.edu.cn}.
}}
\date{}

\begin{document}

\maketitle
\begin{abstract}
In this paper,
we consider a reaction-diffusion-advection SIS epidemic model with saturated incidence rate
and linear source.
We study the uniform bounds of parabolic system and some asymptotic behavior of the
basic reproduction number $\rr_0$,
according to which the threshold-type dynamics are given.
Furthermore,
we show the globally asymptotic stability of the endemic equilibrium in a special case with a large saturated incidence rate.
Meanwhile,
we investigate the effects of diffusion, advection,
saturation and linear source on the asymptotic profiles of the endemic equilibrium. Our study shows that
advection may induce the concentration phenomenon and the disease will be eradicated with small dispersal rate of
infected individuals.
Moreover, our results also reveal that
the linear source can enhance persistence of infectious disease and
large saturation can help speed up the elimination of disease.
\end{abstract}



\section{Introduction}

In this paper,
we consider the following SIS (susceptible-infected-susceptible) reaction-diffusion-advection epidemic model with saturated incidence
\beq
\tag{P}
\label{P}
\left\{\arraycolsep=1.5pt
\begin{array}{lll}
S_t=d_S S_{xx}-q S_x+\Lam(x)-\mu(x)S-\F{\beta(x)SI}{1+mI}+\gam(x)I,\   \ &0<x<L,\ t>0,\\[2mm]
I_t=d_I I_{xx}-q I_x+\F{\beta(x)SI}{1+mI}-\gam(x)I, \   \ &0<x<L,\ t>0,\\[2mm]
d_S S_x-q S=d_I I_x-q I=0,\ &x=0,L,\ t>0,\\[2mm]
S(x,0)=S_0(x)\ge,\not\equiv0,\ I(x,0)=I_0(x)\ge,\not\equiv0, \ &0<x<L,
\end{array}
\right.
\eeq
where $S(x,t)$ and $I(x,t)$ represent the densities of susceptible and infected population at location $x$ and time $t$, respectively;
the positive coefficients $d_S$ and $d_I$, respectively, represent the dispersal rates of susceptible and infected individuals;
the nonnegative constant $q$ denotes the advection speed; the size of the habitat is $L$;
the constant $m>0$ is the saturated incidence rate; $\Lam(x)-\mu(x)S$ indicates that the susceptible individuals are subject to a linear growth,
$\Lam(x)$ and $\mu(x)$ denote the birth and death rates, respectively;
the function $\be(x)$ is the rate of disease transmission,
and $\gam(x)$ is the rate of recovery of infected individuals.
We assume that $\Lam, \mu, \be, \gam$ are positive H\"{o}lder continuous functions for $x\in[0,L]$.
The no-flux boundary conditions mean that no population flux crosses the upstream end($x=0$) and the downstream end($x=L$).

Problem \eqref{P} is a semilinear reaction-diffusion-advection system exhibiting a coupling on the unknowns $S(x,t)$ and $I(x,t)$.
These can be thought of as the descriptions of spatial spread of infectious diseases in the theory of epidemiology.

It is well-known that the spatial spread of diseases in heterogeneous habitats has received considerable attention recently,
since the spatial heterogeneity and individual diffusion can be significant factors in disease dynamics \cite{am1,am2,bc,ddh,dh,h}.
In recent decades,
a lot of works have been
contributed to various reaction-diffusion epidemic models, and to
the effects of spatial heterogeneity and movement of individual on dynamics of infected diseases.

To our knowledge, currently there seems very few works attempting to understand the population dynamics of system \eqref{P}.
In this paper,
we will give a picture of the spread of infectious diseases of model \eqref{P}.
To be more specific,
we first establish some asymptotic behaviors of the basic reproduction number and the relative eigenvalue problem;
see Theorem \ref{properties of R0} and Theorem \ref{lim of lam1} below.
Secondly, we show a global asymptotic stability of  the endemic equilibrium (EE),
which says the EE is globally attractive  regardless of initial values $(S_0,I_0)$;
see Theorem \ref{global stablity of ee}.
We finally focus on the concentration behaviors of endemic equilibrium in six cases: large advection; small diffusion
of the susceptible population; small diffusion of the infected population; large diffusion
of the susceptible population; large diffusion of the infected population; large saturated incidence rate; see Theorem \ref{ee q-infty}-\ref{ee m-infty}.

\subsection{Motivation and related work}

Our work is motivated firstly by \cite{abln}, where they proposed
a classic SIS reaction-diffusion epidemic model:
\beq
\label{classic sis}
\left\{\arraycolsep=1.5pt
\begin{array}{lll}
S_t=d_S\Delta S-\F{\be(x)SI}{S+I}+\gam(x)I,\   \ &\hbox{in} \,\  \Omega\times(0,T),\\[2mm]
I_t=d_I\Delta I+\F{\be(x)SI}{S+I}-\gam(x)I,\   \ &\hbox{in} \,\  \Omega\times(0,T),\\[2mm]
\F{\pa S}{\pa\nu}=\F{\pa I}{\pa\nu}=0,\ &\hbox{on} \ \partial \Omega\times(0,T),\\[2mm]
S(x,0)=S_0(x)\ge,\not\equiv0, I(x,0)=I_0(x)\ge,\not\equiv0 \ &\hbox{in}\,\  \Omega.
\end{array}
\right.
\eeq
It is  shown  in \cite{abln} that
the global dynamics of \eqref{classic sis} can be obtained by the basic reproduction number $\mathcal{R}_0$: the unique disease-free
equilibrium (DFE) is global asymptotically stable if $\mathcal{R}_0<1$ and \eqref{classic sis} admits a unique endemic
equilibrium (EE) if $\mathcal{R}_0>1$.
Moreover,
the authors also investigated the asymptotic behaviors of the endemic equilibrium as the diffusion rate of the susceptible individuals approaches to zero in \cite{abln}.
In \cite{pl},
Peng and Liu discussed the global stability of the endemic equilibrium in some special cases.
In fact, based on \cite{abln},
many subsequence works have been concerned about the dynamics of some related epidemic reaction-diffusion models.
For example, in \cite{p,py},
the authors considered some further asymptotic profiles of the endemic
equilibrium.
Peng and Zhao \cite{pz} considered the same SIS reaction-diffusion model, but the rates of disease transmission and recovery are assumed to be spatially heterogeneous and temporally periodic.
For the reaction-diffusion SIS model with the mass action or
other infection mechanism and the related
epidemic systems, we refer interested readers to \cite{cww,cs1,cs2,dw,gcrfw,lll,lb,lpw,pl,wz,wjl,wz1,zgl} and the references therein.

In some circumstances with external natural forces like
wind, stream flow or water column,
populations may be forced to take passive movement in certain direction.
Concerning this biased dispersal effect,
which is usually modeled by a advection
term,
Cui and Lou \cite{cl} proposed the following SIS epidemic reaction-diffusion-advection model in one-dimensional space domain:
\beq
\label{classic sis advection}
\left\{\arraycolsep=1.5pt
\begin{array}{lll}
S_t=d_S S_{xx}-qS_{x}-\F{\be(x)SI}{S+I}+\gam(x)I,\   \ &0<x<L,t>0,\\[2mm]
I_t=d_I I_{xx}-qI_x+\F{\be(x)SI}{S+I}-\gam(x)I,\   \ &0<x<L,t>0,\\[2mm]
d_S S_x-qS=d_i I_x-qI=0,\ &x=0,L,t>0,\\[2mm]
S(x,0)=S_0(x)\ge,\not\equiv0, I(x,0)=I_0(x)\ge,\not\equiv0, \ &0<x<L.
\end{array}
\right.
\eeq
Compared with the case of no advection,
in \cite{cl},
some different and interesting results on the effects of diffusion and advection rates
on the existence of EE and the
stability of the unique DFE are obtained.
Furthermore, there have been extensive studies
on the diffusion SIS epidemic model in advective environments, see, e.g., \cite{c,cll,cl,clpz,gklz,jwz,kmp,lz,zc}.

Notice that a very important feature of model \eqref{classic sis} and \eqref{classic sis advection} is that they have the conservation law,
i.e.,
the total number of individuals in their habitats is conserved provided that the initial total population is fixed;
see \cite{abln,cl}.
Therefore, it is natural to consider the scenario that the susceptible individuals are allowed to have birth and death.
Recently, Chen and Cui \cite{cc} studied
the following reaction-diffusion-advection SIS epidemic model with linear source,
which has the vary total individual:
\beq
\label{chen cui}
\left\{\arraycolsep=1.5pt
\begin{array}{lll}
S_t=d_S S_{xx}-q S_x+\Lam(x)-\mu(x)S-\F{\be(x)SI}{S+I}+\gam(x)I,\   \ &0<x<L,\ t>0,\\[2mm]
I_t=d_I I_{xx}-q I_x+\F{\be(x)SI}{S+I}-\gam(x)I, \   \ &0<x<L,\ t>0,\\[2mm]
d_S S_x-q S=d_I I_x-q I=0,\ &x=0,L,\ t>0,\\[2mm]
S(x,0)=S_0(x)\ge,\not\equiv0,\ I(x,0)=I_0(x)\ge,\not\equiv0, \ &0<x<L.
\end{array}
\right.
\eeq
Here, $S,I,d_S,d_I,q,\Lam,\mu,\be,\gam$ have the same meaning as in \eqref{P}.

In fact,
our work is also motivated by \cite{cc},
where they derived the persistence and extinction
of the infectious disease in terms of the
basic reproduction number.
Moreover, the global
asymptotic stability of the endemic equilibrium in a special case
and the asymptotic profiles of the endemic equilibrium with respect to linear
source, advection, diffusion are also discussed.

We end this subsection by reviewing an open problem proposed in \cite{cc}.
When the infection mechanism is governed by the well-known mass action mechanism
or saturated incidence mechanism, instead of
standard incidence infection one,
what is the dynamics of the corresponding systems?
The purpose of this paper is to pursue this subject,
i.e.,
we shall study model \eqref{P}.

\subsection{Main results}

To investigate the dynamics of system \eqref{P},
we shall introduce non-negative equilibrium solutions of \eqref{P},
that is, the non-negative solutions of the following system:
\beq
\label{equilibrium}
\left\{\arraycolsep=1.5pt
\begin{array}{ll}
d_S S_{xx}-q S_x+\Lam(x)-\mu(x)S-\F{\be(x)SI}{1+mI}+\gam(x)I=0,\   \ &0<x<L,\\[2mm]
d_I I_{xx}-q I_x+\F{\be(x)SI}{1+mI}-\gam(x)I=0, \   \ &0<x<L,\\[2mm]
d_S S_x-q S=d_I I_x-q I=0,\ &x=0,L.
\end{array}
\right.
\eeq
Clearly, we are interested only in solutions $(S(x),I(x))$ satisfying $S(x)\ge0$ and $I(x)\ge0$ on $[0,L]$.
We denote an EE by $(\tilde{S},\tilde{I})$.
It is evident that \eqref{equilibrium}
admits an unique disease-free equilibrium (DFE) $(\hat{\tilde{S}},0)$,
where $\hat{\tilde{S}}$ is the unique positive solution to
\beq
\label{hattilde S}
\left\{\arraycolsep=1.5pt
\begin{array}{ll}
d_S S_{xx}-q S_x+\Lam(x)-\mu(x)S=0,\   \ &0<x<L,\\[2mm]
d_S S_x-q S=0,\ &x=0,L.
\end{array}
\right.
\eeq
From a similar manner in literature \cite{cc,cl,wz},
we can derive the basic reproduction number $\mathcal{R}_0$ of the system
\eqref{P} as follows:
\beq
\label{R0}
\mathcal{R}_0(d_I,d_S,q)=
\sup\limits_{\varphi\in H^1((0,L)),\varphi\neq0}
\F{\i_0^L\s(x)\be(x)e^{\f{q}{d_I}x}\varphi^2\dx}{d_I\i_0^Le^{\f{q}{d_I}x}\varphi_x^2\dx+\i_0^L\gam(x)e^{\f{q}{d_I}x}\varphi^2\dx}.
\eeq

We are now in a position to state our  main results.
Note by the standard theory for semilinear parabolic systems
that
system \eqref{P} admits a unique classical solution ($S(x,t),I(x,t)$) \cite{h0}.
Moreover, it follows from the maximum principle \cite{pw,w} that
both $S(x,t)$ and $I(x,t)$ are positive for $(x,t)\in[0,L]\times(0,+\infty)$.
Throughout our paper,
we use $\|\cdot\|_p$ to  denote  the standard norm of
$L^{p}((0,L))$, and denote
\beq
f^*=\max\limits_{x\in[0,L]}f(x),\ f_*=\min\limits_{x\in[0,L]}f(x).
\eeq

\begin{thm}
\label{properties of R0}
For any $d_S>0,d_I>0,q>0$,
we have the following statements:

(1)
\beq
\begin{cases}
\mathcal{R}_0<1\Leftrightarrow\lam_1(d_I,q,\gam-\s\be)>0,\\
\mathcal{R}_0=1\Leftrightarrow\lam_1(d_I,q,\gam-\s\be)=0,\\
\mathcal{R}_0>1\Leftrightarrow\lam_1(d_I,q,\gam-\s\be)<0,
\end{cases}
\eeq
where $\lam_1(d_I,q,\gam-\s\be)$ is the principal eigenvalue of the eigenvalue problem:
\beq
\label{lam1(d q ga-s)}
\left\{\arraycolsep=1.5pt
\begin{array}{ll}
-d_I\psi_{xx}+q\psi_x+(\gam(x)-\s(x)\be(x))\psi=\lam\psi,\   \ &0<x<L,\\[2mm]
d_I\psi_x-q\psi=0,\ &x=0,L,
\end{array}
\right.
\eeq
and can be
characterized as
\beq
\label{lam1 variational}
\lam_1(d_I,q,\gam-\s\be)
=\inf\limits_{0\neq\eta\in H^1_0(0,L)}
\F{d_I\io e^{\f{qx}{d_I}}(\eta_x)^2\dx+\io(\gam(x)-\be(x)\s(x))e^{\f{qx}{d_I}}\eta^2(x)\dx}
{\io e^{\f{qx}{d_I}}\eta^2(x)\dx};
\eeq

(2) Given any $d_S>0,q>0$,
$\mathcal{R}_0\rightarrow\F{\s(L)\be(L)}{\gam(L)}$ as $d_I\rightarrow0$;

(3) Given any $d_S>0,d_I>0$,
$\mathcal{R}_0\rightarrow+\infty$ as $q\rightarrow+\infty$;

(4) Given any $d_I>0,q>0$,
$\rr_0\rightarrow\rr_0^*:=\F{\|\Lam\|_1}{\|\mu\|_1}\sup\limits_{\varphi\in H^1((0,L)),\varphi\neq0}\F{\i_0^L\be(x)e^{\f{q}{d_I}x}\varphi^2\dx}{d_I\i_0^Le^{\f{q}{d_I}x}\varphi_x^2\dx+\i_0^L\gam(x)e^{\f{q}{d_I}x}\varphi^2\dx}$ as $d_S\rightarrow+\infty$.
\end{thm}

\begin{rem}
Theorem \ref{properties of R0} (3) implies that the basic reproduction number $\rr_0$ is always larger than $1$ provided that the advection is sufficiently large.
\end{rem}

\begin{rem}
Indeed, if $\Lam(x)\equiv\mu(x)$,
then $\rr_0^*$ is equal to the basic reproduction number of the classical SIS model \eqref{classic sis advection} with standard incidence and without linear source.
\end{rem}

The principal eigenvalue $\lam_1(d_I,q,\gam-\s\be)$ of \eqref{lam1(d q ga-s)} has the following limiting profile.

\begin{thm}
\label{lim of lam1}
Given $d_I,q>0$,
if $\mu'(x)\le0,\Lam'(x)\ge0$,
then there holds
$\lim\limits_{d_S\rightarrow0}\lam_1(d_I,q,\gam-\s\be)=\bar{\lam}$,
where $\bar{\lam}$ is the principal eigenvalue of
\beq
\label{lim eigen problem}
\left\{\arraycolsep=1.5pt
\begin{array}{ll}
-d_I\bar{\psi}_{xx}+q\bar{\psi}_x+(\gam(x)-\be\s^{\infty})\bar{\psi}=\bar{\lam}\bar{\psi},\   \ &0<x<L,\\[2mm]
d_I\bar{\psi}_x(0)-q\bar{\psi}(0)=0,\ &\\[2mm]
d_I\bar{\psi}_x(L)-q\bar{\psi}(L)=\be(L)\bar{\psi}(L)\mu(L)N_S.&
\end{array}
\right.
\eeq
Here
$\s^{\infty}(x)$ is the solution of
\beq
\label{s infty}
\left\{\arraycolsep=1.5pt
\begin{array}{ll}
-q\s^{\infty}_x+\Lam(x)-\mu(x)\s^{\infty}=0,\   \ &0<x<L,\\[2mm]
\s^{\infty}(0)=0,\ &
\end{array}
\right.
\eeq
and
\beq
\label{ns}
N_S:=\F{\|\Lam\|_1-\Io\mu(x)\s^{\infty}(x)\dx}{\mu(L)}.
\eeq
\end{thm}

% \begin{rem}
% When the minimal solution $(w_{\lam,\mu},z_{\lam,\mu})$ is the only solution to \eqref{E}, except for the case in Theorem \ref{Global existence vs quenching} $(c)$, the behavior of solution  $(u(x,t),v(x,t))$ to \eqref{P} with initial data above $(w_{\lam,\mu},z_{\lam,\mu})$ is unknown.
% \end{rem}

Next, before we give a globally asymptotically
stability of a special EE,
we need the following assumption:
\beq
\label{assumption}
\begin{cases}
\Lam(x)\be(x)-\gam(x)\mu(x)>0,\\
\F{\Lam(x)}{\mu(x)e^{\f{qx}{d_S}}}=\kap\ \text{for some constant}\ \kap\in(0,+\infty),\\
\F{\Lam(x)\be(x)-\gam(x)\mu(x)}{\gam(x)\mu(x)e^{\f{qx}{d_I}}}=r\ \text{for some constant}\ r\in(0,+\infty).
\end{cases}
\eeq
Then it is easy to see that
\beq
(\tilde{S},\tilde{I})=\big(\F{\Lam(x)}{\mu(x)},\F{\Lam(x)\be(x)-\gam(x)\mu(x)}{\gam(x)\mu(x)}\big)
\eeq
is a EE of \eqref{P}.

\begin{thm}
\label{global stablity of ee} Suppose that \eqref{assumption} hold.
Then there exists some $M>0$,
such that for any $m>M$,
the EE is globally attractive.
\end{thm}

% \begin{rem} \label{7103}
% Assume further that $(u_0,v_0)\leq(w_{\lam,\mu},z_{\lam,\mu})$
% and $(u_0,v_0)$ is a subsolution of \eqref{E}, or that $(w_{\lam,\mu},z_{\lam,\mu})\leq(u_0,v_0)\le\not\equiv
% (w_1, z_1)$ (if $(w_1,z_1)$ exists) and  $(u_0,v_0)$ is a   supersolution of \eqref{E}. Then the exponential convergence further holds in $H^1$ norm. More precisely, there holds
%  \beq
%   \|u(x,t)-w_{\lam,\mu}(x)\|_{H^1}+\|v(x,t)-z_{\lam,\mu}(x)\|_{H^1}\leq C \exp\bigg({-\min\bigg\{\frac{\lambda_1}{4},\frac{\nu_1}{4},\frac{1}{2}\bigg\}t}\bigg),\quad \ t>T_0
%   \eeq
%   for some $T_0>0$. Here, $C$ is a constant depending at most on $u_0,v_0,w_{\lam,\mu},z_{\lam,\mu},w_1,z_1,\varphi_1,\psi_1;$  $\lambda_1$ and $\nu_1$ are defined as in Theorem \ref{Convergence rate}.
% \end{rem}

The following results demonstrate some concentration behaviors of EE.

\begin{thm}
\label{ee q-infty}
Assume that $\rr_0>1$.
Let $(S(x),I(x))$ be any positive solution of \eqref{equilibrium}.
Given $d_S,d_I,m>0$,
as $q\rightarrow+\infty$,
the following statements hold:

(1)
$\lim\limits_{q\rightarrow+\infty}(S(x),I(x))=(0,0)$ locally uniformly in $[0,L)$;

(2)
$\lim\limits_{q\rightarrow+\infty}(\F{1}{q}S(L-\F{y}{q}),\F{1}{q}I(L-\F{y}{q}))=(\F{\|\Lam\|_1}{d_S\mu(L)}e^{-\f{y}{d_S}},\F{\be(L)\|\Lam\|_1}{m\gam(L)d_I\mu(L)}e^{-\f{y}{d_I}})$,
locally uniformly for $y\in[0,+\infty)$.
\end{thm}

Next, before giving the asymptotic profile of EE as $d_S\rightarrow0$,
we need introduce a eigenvalue problem as follows:
\beq
\label{aux eigen}
\left\{\arraycolsep=1.5pt
\begin{array}{ll}
-d_I\xi_{xx}+q\xi_x+(\gam(x)-\s^{\infty}(x)\be(x))\xi=\tau\xi,\   \ &0<x<L,\\[2mm]
d_I\xi_x(0)-q\xi(0)=0,\ &\\[2mm]
d_I\xi_x(L)-q\xi(L)=n\be(L)\mu(L)\xi(L),&
\end{array}
\right.
\eeq
where $n>0$,
$\s^{\infty}$ is the solution of
\eqref{s infty}.
Indeed, the principal eigenvalue of \eqref{aux eigen},
denoted by $\tau_1(n)$,
can be
characterized as
\beq
\tau_1(n)=\inf\limits_{0\neq\omega\in H^1_0(0,L)}
\F{-n\be(L)\mu(L)e^{\f{qL}{d_I}}\omega^2(L)+d_I\io e^{\f{qx}{d_I}}(\omega_x)^2\dx+\io(\gam(x)-\be(x)\s^{\infty}(x))e^{\f{qx}{d_I}}\omega^2(x)\dx}{\io e^{\f{qx}{d_I}}\omega^2(x)\dx}.
\eeq
Moreover, $\tau(n)$ is smooth and decreasing in $n$.
%Similar to the eighth line on Page 6 in \cite{clpz},
Using the variational characterization above,
we get that
\beq
\tau_1(n)\le
\F{-n\be(L)\mu(L)e^{\f{qL}{d_I}}+\io(\gam(x)-\be(x)\s^{\infty}(x))e^{\f{qx}{d_I}}\dx}
{\io e^{\f{qx}{d_I}}\dx}\le
\F{-n\be(L)\mu(L)e^{\f{qL}{d_I}}+\io\gam(x)e^{\f{qx}{d_I}}\dx}
{\io e^{\f{qx}{d_I}}\dx},
\eeq
which implies that $\lim\limits_{n\rightarrow+\infty}\tau_1(n)=-\infty$.
Note that
\beq
\bea{l}
\tau_1(0)=\inf\limits_{0\neq\omega\in H^1_0(0,L)}
\F{d_I\io e^{\f{qx}{d_I}}(\omega_x)^2\dx+\io(\gam(x)-\be(x)\s^{\infty}(x))e^{\f{qx}{d_I}}\omega^2(x)\dx}{\io e^{\f{qx}{d_I}}\omega^2(x)\dx}\\
\begin{cases}
\ge\inf\limits_{0\neq\omega\in H^1_0(0,L)}
\F{\io(\gam(x)-\be(x)\s^{\infty}(x))e^{\f{qx}{d_I}}\omega^2(x)\dx}{\io e^{\f{qx}{d_I}}\omega^2(x)\dx}>0,\ \text{if}\ \F{\gam(x)}{\be(x)}>\s^{\infty}(x)\ \text{for all}\ x\in[0,L],\\
\le\F{\io(\gam(x)-\be(x)\s^{\infty}(x))e^{\f{qx}{d_I}}\dx}{\io e^{\f{qx}{d_I}}\dx}\le0,\ \text{if}\ \F{\gam(x)}{\be(x)}\le\s^{\infty}(x)\ \text{for all}\ x\in[0,L].
\end{cases}
\eea
\eeq
Therefore,
one can have that there exists a unique $N_0\in(0,+\infty)$ such that
\beq
\tau_1(n)
\begin{cases}
>0,\ \text{when}\ 0<n<N_0,\\
=0,\ \text{when}\ n=N_0,\\
<0,\ \text{when}\ n>N_0,
\end{cases}
\eeq
provided that $\F{\gam(x)}{\be(x)}>\s^{\infty}(x)$.
On the other hand,
$\tau_1(n)<0$ for any $n>0$,
provided that $\F{\gam(x)}{\be(x)}\le\s^{\infty}(x)$.

\begin{rem}
In fact, one can easily see that $\bar{\lam}=\tau_1(N_S)$.
\end{rem}

\begin{rem}
Similar to Theorem \ref{properties of R0} (1),
it is easy to see that
\beq
\begin{cases}
\lim\limits_{d_S\rightarrow0}\rr_0<1\Leftrightarrow\tau_1(N_S)>0,\\
\lim\limits_{d_S\rightarrow0}\rr_0=1\Leftrightarrow\tau_1(N_S)=0,\\
\lim\limits_{d_S\rightarrow0}\rr_0>1\Leftrightarrow\tau_1(N_S)<0.
\end{cases}
\eeq
\end{rem}

\begin{thm}
\label{ee ds-0}
Suppose that $d_I,q,m>0$.
Assume that
\beq
\label{assum}
\text{either}\ N_S>N_0\ \text{and}\ \F{\gam(x)}{\be(x)}>\s^{\infty}(x)\ \text{or}\ \F{\gam(x)}{\be(x)}\le\s^{\infty}(x).
\eeq
Let $(S(x),I(x))$ be any EE of \eqref{equilibrium}.
Then as $d_S\rightarrow0$,
if $\mu'(x)\le0,\Lam'(x)\ge0$,
we have

(1) $d_S S(L-d_S y)
\rightarrow\F{q\|\Lam\|_1}{\mu(L)} e^{-qy}$
locally uniformly in $[0,+\infty)$;

(2)
$(S(x),I(x))\rightarrow(S^{\infty}(x)+N_S\del_{L}(x),I^{\infty}(x))$ weakly in $L^1((0,L))\times L^1((0,L))$,
where $\del_L(x)$ stands for the Dirac delta distribution at $L$,
$N_S$ is as in \eqref{ns},
$(S^{\infty}(x),I^{\infty}(x))$ is a positive solution of
\beq
\label{sinfty iinfty}
\left\{\arraycolsep=1.5pt
\begin{array}{ll}
-q S^{\infty}_x+\Lam(x)-\mu(x)S^{\infty}-\F{\be(x)S^{\infty}I^{\infty}}{1+mI^{\infty}}+\gam(x)I^{\infty}=0,\   \ &0<x<L,\\[2mm]
d_I I^{\infty}_{xx}-q I_x+\F{\be(x) S^{\infty}I^{\infty}}{1+mI^{\infty}}-\gam(x)I^{\infty}=0, \   \ &0<x<L,\\[2mm]
S^{\infty}(0)=0,\ &,\\[2mm]
d_I I^{\infty}_x(0)-qI^{\infty}(0)=0,&\\[2mm]
d_I I_x^{\infty}(L)-qI^{\infty}(L)=N_S\F{\be(L)I^{\infty}(L)}{1+mI^{\infty}(L)}.&
\end{array}
\right.
\eeq
\end{thm}

\begin{thm}
\label{ee di-0}
Suppose that $\rr_0>1$.
Let $(S(x),I(x))$ be any positive solution of \eqref{equilibrium}.
Then for any $d_S,q>0$,

(1) $\lim\limits_{d_I\rightarrow0}S(x)=\s(x)$ uniformly on $[0,L]$;

(2)
$\lim\limits_{d_I\rightarrow0}I(x)=0$ locally uniformly in $[0,L)$,
and $\lim\limits_{d_I\rightarrow0}\Io I(x)\dx=0$.
\end{thm}

\begin{thm}
\label{ee ds-infty}
Suppose that $\rr_0^*>1$.
Let $(S(x),I(x))$ be any EE of \eqref{equilibrium}.
Then as $d_S\rightarrow+\infty$, we have
\beq
\lim\limits_{d_S\rightarrow+\infty}(S(x),I(x))=(\F{\|\Lam\|_1}{\|\mu\|_1},I^{\infty}(x))\ \text{uniformly on}\ [0,L].
\eeq
Here $I^{\infty}$ is the unique positive
solution of
\beq
\left\{\arraycolsep=1.5pt
\begin{array}{ll}
d_I I^{\infty}_{xx}-q I^{\infty}_x+\beta(x)\F{S^{\infty}I^{\infty}}{1+mI^{\infty}}-\gam(x)I^{\infty}=0, \   \ &0<x<L,\\[2mm]
d_I I^{\infty}_x-q I^{\infty}=0,\ &x=0,L.
\end{array}
\right.
\eeq
\end{thm}

\begin{thm}
\label{ee di-infty}
Assume that $\rr_0>1$, $\Io\be(x)\s(x)\dx>\Io\gam(x)\dx$.
Let $(S(x),I(x))$
be any positive solution of \eqref{equilibrium}.
Fix $d_S,q,m>0$,
then
\beq
\lim\limits_{d_I\rightarrow+\infty}(S(x),I(x))=(S_{\infty}(x),I_{\infty})\ \text{uniformly on}\ [0,L],
\eeq
where $I_{\infty}$ is a positive constant and $S_{\infty}>0$ solves
\beq
\left\{\arraycolsep=1.5pt
\begin{array}{ll}
d_S (S_{\infty})_{xx}-q (S_{\infty})_x+\Lam(x)-\mu(x)S_{\infty}-\F{\be(x)S_{\infty}I_{\infty}}{1+mI_{\infty}}+\gam(x)I_{\infty}=0, \   \ &0<x<L,\\[2mm]
d_S (S_{\infty})_x-q S_{\infty}=0,\ &x=0,L,\\[2mm]
\Io\gam(x)\dx=\Io\F{\be(x)S_{\infty}(x)}{1+mI_{\infty}}\dx.
\end{array}
\right.
\eeq
\end{thm}

\begin{thm}
\label{ee m-infty}
Suppose that $\rr_0>1$.
Let $(S(x),I(x))$ be any positive solution of \eqref{equilibrium}.
Then for any $d_S,d_I,q>0$,
$\lim\limits_{m\rightarrow+\infty}(S(x),I(x))=(\s(x),0)$ uniformly on $[0,L]$,
where $\s$ is the unique solution of \eqref{hattilde S}.
Moreover, $\lim\limits_{m\rightarrow+\infty}mI=\th^*(x)$,
where $\th^*$ is the unique solution of
\beq
\label{th*}
\left\{\arraycolsep=1.5pt
\begin{array}{ll}
d_I \th^*_{xx}-q \th^*_x+(\F{\be \s}{1+\th^*}-\gam(x))\th^*=0, \   \ &0<x<L,\\[2mm]
d_I \th^*_x-q \th^*=0,\ &x=0,L.
\end{array}
\right.
\eeq
\end{thm}

We  remark that one of the main difference between our paper and \cite{cc}
is that the term $\F{S}{1+mI}$ may not bounded.
$\F{S}{S+I}$,
however,
is bounded,
which plays an important role in \cite{cc}.
Hence,
in this paper,
we have to use some different ways to consider our model.

This paper is organized as follows.
In Section \ref{Uniform boundedness},
we provide some uniform boundedness of the solution of \eqref{P}.
In Section \ref{Extinction and persistence},
we display the threshold-type dynamics in terms of the
basic reproduction number $\rr_0$
and verify some properties of $\rr_0$.
In Section \ref{Global attractivity of the EE},
we prove Theorem \ref{global stablity of ee}.
Some concentrations of EE are given in Section \ref{Asymptotic profiles of the EE}.
At last,
Section \ref{discussion} is devoted to a brief discussion about our results.

\section{Uniform boundedness}
\label{Uniform boundedness}

In this section,
we provide the uniform boundedness of solutions to \eqref{P}.

\begin{lem}
\label{Uniform bounds}
Assume that $d_S=d_I$.
For any solution $(S,I)$ of \eqref{P},
there holds
\beq
S(x,t)+I(x,t)\le
\max\Big\{\Big(\F{\Lam}{\sig e^{\f{qx}{d}}}\Big)^*e^{\f{qx}{d}},
\Big(\f{S_0+(1+\eps_0)I_0}{e^{\f{qx}{d}}}\Big)^*e^{\f{qx}{d}}\Big\}.
\eeq
Here $\eps_0$ is any given positive constant satisfying $0<\eps_0<\F{m\mu_*}{\be^*}$
and $\sig=\min\{\F{\eps_0\gam_*}{1+\eps_0},\mu_*-\F{\eps_0\be^*}{m}\}$.
\end{lem}

\proof
We let $d=d_S=d_I$
and set $H(x,t)=S(x,t)+(1+\eps_0)I(x,t)$,
where $\eps_0$ is any given positive constant satisfying $0<\eps_0<\F{m\mu_*}{\be^*}$.
By direct calculations,
we reach that
\beq
\bea{l}
\F{\pa H}{\pa t}=d H_{xx}-q H_x+\Lam(x)-\mu(x)S+\eps_0\be(x)\F{SI}{1+mI}-\eps_0\gam(x)I\\
\le d H_{xx}-q H_x+\Lam(x)-\mu_*S+\F{\eps_0\be^*}{m}S-\eps_0\gam_*I\\
\le d H_{xx}-q H_x+\Lam(x)-\sig H.
\eea
\eeq
%where $\sig=\min\{\F{\eps_0\gam_*}{1+\eps_0},\mu_*-\F{\eps_0\be^*}{m}\}$.
Since $\max\Big\{(\F{\Lam}{\sig e^{\f{qx}{d}}})^*e^{\f{qx}{d}},(\f{S_0+(1+\eps_0)I_0}{e^{\f{qx}{d}}})^*e^{\f{qx}{d}}\Big\}$ and $H$ is a pair of upper and lower solutions to the
initial-boundary value problem:
\beq
\left\{\arraycolsep=1.5pt
\begin{array}{lll}
u_t=d u_{xx}-q u_x+\Lam(x)-\sig u,\   \ &0<x<L,\ t>0,\\[2mm]
d u_x-q u=0,\ &x=0,L,\ t>0,\\[2mm]
u(x,0)=S_0(x)+(1+\eps_0)I_0(x)\ge,\not\equiv0, \ &0<x<L,
\end{array}
\right.
\eeq
it can be deduced by the comparison principle that for any $x\in[0,L]$ and $t>0$,
\beq
H(x,t)\le\max
\Big\{(\F{\Lam}{\sig e^{\f{qx}{d}}})^*e^{\f{qx}{d}},(\f{S_0+(1+\eps_0)I_0}{e^{\f{qx}{d}}})^*e^{\f{qx}{d}}\Big\}.
\eeq
\proofend

\begin{lem}
\label{Ultimate uniform bounds}
For any solution $(S,I)$ of \eqref{P},
there holds for every $t>0$,
\beq
\|S(\cdot,t)\|_{\infty}+\|I(\cdot,t)\|_{\infty}\le M_1,
\eeq
for some $M_1>0$ depending on
initial value.
Moreover, there exists some positive constant $M_2$ independent of initial data such that
\beq
\label{uniform bound}
\|S(\cdot,t)\|_{\infty}+\|I(\cdot,t)\|_{\infty}\le M_1,\ \forall t\ge T,
\eeq
for some large $T$.
\end{lem}

\proof
Let $\eps_0$ be a given positive constant fulfilling $0<\eps_0<\F{\mu_*m}{\be^*}$,
and set
\beq
H(t)=\Io(S(x,t)+(1+\eps_0)I(x,t))\dx.
\eeq
It then follows that
\beq
\begin{array}{ll}
H'(t)&=\Io\Lam(x)\dx-\Io\mu(x)S(x)\dx
+\eps_0\Io\F{\be(x)SI}{1+mI}\dx-\eps_0\Io\gam(x)I\dx\\
&\le\|\Lam\|_1-\mu_*\Io S\dx+\F{\eps_0\be^*}{m}\Io S\dx-\eps_0\gam_*\Io I\dx\\
&\le\|\Lam\|_1-\sig H(t),
\end{array}
\eeq
where
$\sig=\min\{\F{\eps_0\gam_*}{1+\eps_0},\mu_*-\F{\eps_0\be^*}{m}\}$.
Then by Gronwall's inequality,
we get that
\beq
H(t)\le H(0)e^{-\sig t}+\F{\|\Lam\|_1}{\sig}(1-e^{-\sig t}),\ \text{for any}\ t\ge0.
\eeq
This implies that any solution $(S(x,t),I(x,t))$ satisfies $L^1$ bound uniformly for all $t\in[0,+\infty)$.
By means of a similar Moser-type iterative method in \cite{a,pz},
we have that both $\|S(\cdot,t)\|_{\infty}$ and $\|I(\cdot,t)\|_{\infty}$ are uniformly bounded for all $t\in[0,+\infty)$.
Making use the method of \cite[Lemma 2.1]{dp},
we can derive \eqref{uniform bound}
by noting that $\limsup\limits_{t\rightarrow+\infty}H(t)\le\F{\|\Lam\|_1}{\sig}$.
The proof is completed.
% \beq
% \left\{
% \bea{l}
% f'(\eta_z)-f'(\xi_z)<0,
% g'(\eta_z)-g'(\xi_z)<0,\\
% w_3\ge w_2\ge w_1,z_3\ge z_2\ge z_1.
% \eea
% \right.
% \eeq
\proofend

\section{Extinction and persistence/the properties of $\rr_0$} \label{Extinction and persistence}

This section is devoted to investigating the properties of the basic reproduction
number $\mathcal{R}_0$ defined in \eqref{R0}
and showing the properties of the extinction and persistence in terms of $\mathcal{R}_0$.
We first give some properties of the basic reproduction number $\mathcal{R}_0$ and relative principal eigenvalue $\lam_1(d_I,q,\gam-\s\be)$,
i.e.,
we shall show Theorem \ref{properties of R0} and Theorem \ref{lim of lam1}.

\one
Indeed, the statements (1) and (2) can be proved similarly in \cite[Lemma 2.2 and Theorem 1.2(i)]{cl}.
Now we will verify statements (3) and (4).

By the definition of the basic reproduction number $\rr_0$ in \eqref{R0},
it follows that
there exists a positive function $\varphi\in C^2([0,L])$ such that
\beq
\label{equation of r0}
\left\{\arraycolsep=1.5pt
\begin{array}{ll}
-d_I\varphi_{xx}+q\varphi_x+\gam(x)\varphi=\F{1}{\rr_0}\be\s\varphi,\   \ &0<x<L,\\[2mm]
d_I\varphi_x-q\varphi=0,\ &x=0,L;
\end{array}
\right.
\eeq
and $\rr_0$ also admits the following
variational characterization:
\beq
\label{r0}
\rr_0(d_I,d_S,q)=
\sup\limits_{\phi\in H^1((0,L)),\phi\neq0}
\F{\i_0^L\s(x)\be(x)\phi^2\dx}{d_I\i_0^L|\phi_x-\f{q}{2d_I}\phi|^2\dx+\i_0^L\gam(x)\phi^2\dx}.
\eeq
Let $\phi_q(x)=\sqrt{\F{q}{d_I(e^{\f{qL}{d_I}}-1)}}e^{\f{qx}{2d_I}}$.
Using \eqref{r0},
we
arrive at
\beq
\begin{array}{l}
\rr_0\ge\F{\i_0^L\s(x)\be(x)\phi_q^2\dx}{d_I\i_0^L|(\phi_q)_x-\f{q}{2d_I}\phi_q|^2\dx+\i_0^L\gam(x)\phi_q^2\dx},\\
=\F{\i_0^L\s(x)\be(x)\phi_q^2\dx}{\i_0^L\gam(x)\phi_q^2\dx}\\
\ge\F{\be_*\io\underline{S}(x)\phi_q^2\dx}{\gam^*},
\end{array}
\eeq
where $\underline{S}(x)$ is the unique positive solution to
\beq
\left\{\arraycolsep=1.5pt
\begin{array}{lll}
d_S S_{xx}-q S_x+\Lam_*-\mu^*S=0,\   \ &0<x<L,\\[2mm]
d_S S_x-q S=0,\ &x=0,L.
\end{array}
\right.
\eeq
After some direct calculations,
it follows that
$\lim\limits_{q\rightarrow+\infty}\F{\be_*\io\underline{S}(x)\phi_q^2\dx}{\gam^*}=+\infty$.
Hence the statement (3) holds.

At last, we shall prove statement (4).
Integrating firstly the equation of $\s$ in \eqref{hattilde S} over $(0,L)$,
there holds
\beq
\label{5}
\Io\mu(x)\s(x)\dx=\Io\Lam(x)\dx,
\eeq
which implies that
\beq
\label{4}
\F{\|\Lam\|_1}{\mu^*}
\le\Io\s(x)\dx\le\F{\|\Lam\|_1}{\mu_*}.
\eeq
Then integrating \eqref{hattilde S} from $0$ to $x$,
together with \eqref{5},
yields that
\beq
d_S\s_x-q\s=\I_0^x(-\Lam(x)+\mu(x)\s(x))\dx
\le\Io\mu(x)\s(x)\dx=\|\Lam\|_1.
\eeq
We multiply both sides of the above inequality by
$e^{-\f{qx}{d_S}}$
and integrate over $(x,L)$ to observe that
\beq
\label{8}
e^{\f{q(x-L)}{d_S}}\s(L)
\le\s(x)+\F{\|\Lam\|_1}{q}(1-e^{\f{q}{d_S}(x-L)}).
\eeq
By further integrating from $0$ to $L$ and direct calculations,
we get that
\beq
\label{s(L)<=}
\bea{l}
\s(L)\le\F{q}{d_S(1-e^{-\f{qL}{d_S}})}\big[\Io\s(x)\dx+\F{\|\Lam\|_1}{q}(L-\F{d_S}{q}(1-e^{-\f{qL}{d_S}}))\big]\\
\le\F{1}{d_S(1-e^{-\f{qL}{d_S}})}\big[\F{q\|\Lam\|_1}{\mu_*}+\|\Lam\|_1(L-\F{d_S}{q}(1-e^{-\f{qL}{d_S}}))\big].
\eea
\eeq

Secondly, integrating \eqref{hattilde S} from $x$ to $L$,
together with \eqref{5},
yields that
\beq
-d_S\s_x+q\s=\I_x^L(-\Lam(x)+\mu(x)\s(x))\dx
\le\Io\mu(x)\s(x)\dx=\|\Lam\|_1.
\eeq
Similarly,
we multiply both sides of the above inequality by
$e^{-\f{qx}{d_S}}$
and integrate over $(x,L)$ to get that
\beq
d_S\s e^{-\f{qx}{d_S}}-d_S\s(L)e^{-\f{qL}{d_S}}
\le\F{d_S\|\Lam\|_1}{q}(e^{-\f{qx}{d_S}}-e^{-\f{qL}{d_S}}),
\eeq
which combined with \eqref{s(L)<=} indicates that
\beq
\label{9}
\bea{l}
\s(x)\le\F{\|\Lam\|_1}{q}
+\F{e^{-\f{q(L-x)}{d_S}}}{d_S(1-e^{-\f{qL}{d_S}})}\big[\F{q\|\Lam\|_1}{\mu_*}+\|\Lam\|_1L\big].
\eea
\eeq
This implies that $\s(x)$ is bounded uniformly on $[0,L]$ for any $d_S>1$.
Integrating \eqref{hattilde S} from $x$ to $L$ again,
one can derive that
\beq
\label{aa}
d_S\s_x-q\s=\I_x^L(\Lam(x)-\mu(x)\s(x))\dx
\le\|\Lam\|_1,
\eeq
which indicates that $\s_x$ is uniformly bounded on $[0,L]$ for $d_S>1$.
The equation \eqref{hattilde S} then yields that $\s_{xx}$ is also uniformly bounded on $[0,L]$ and continuous for $d_S>1$.
By means of the Sobolev embedding theorem,
by passing to the subsequences of $d_S\rightarrow+\infty$ if necessary,
we have that
$\s(x)\rightarrow S^{\infty}(x)$ in $C^1([0,L])$.
Dividing both side of \eqref{aa} by $d_S$ and sending $d_S\rightarrow+\infty$ give rise to $S^{\infty}_x\equiv 0$.
So $S^{\infty}(x)$ is a constant with this constant being $\F{\|\Lam\|_1}{\|\mu\|_1}$ by using \eqref{5}.
Therefore $\lim\limits_{d_S\rightarrow+\infty}\s(x)=\F{\|\Lam\|_1}{\|\mu\|_1}$.
Then
\beq
\lim\limits_{d_S\rightarrow+\infty}\rr_0=\rr_0^*:=\F{\|\Lam\|_1}{\|\mu\|_1}\sup\limits_{\varphi\in H^1((0,L)),\varphi\neq0}\F{\i_0^L\be(x)e^{\f{q}{d_I}x}\varphi^2\dx}{d_I\i_0^Le^{\f{q}{d_I}x}\varphi_x^2\dx+\i_0^L\gam(x)e^{\f{q}{d_I}x}\varphi^2\dx}.
\eeq
\proofend

\oneplus
We aim to show this theorem in several steps.

\emph{Step 1.}
Given $d_I,q>0$,
$\liminf\limits_{d_S\rightarrow0}\lam_1(d_I,q,\gam-\s\be)>-\infty$.

By a transformation $\omega=e^{-\f{qx}{2d_S}}\psi$,
we can get the following equivalent form of the eigenvalue problem \eqref{lam1(d q ga-s)}:
\beq
\left\{\arraycolsep=1.5pt
\begin{array}{ll}
-d_I\omega_{xx}+(\F{q^2}{4d_I}+\gam(x)-\s(x)\be(x))\omega=\lam_1\omega,\   \ &0<x<L,\\[2mm]
d_I\omega_x=\F{q}{2}\omega,\ &x=0,L.
\end{array}
\right.
\eeq
Define now a new eigenvalue problem:
\beq
\label{new eigenvalue problem}
\left\{\arraycolsep=1.5pt
\begin{array}{ll}
-d_I\omega_{xx}+(\F{q^2}{4d_I}+\gam_*-\be^*\s(x))\omega=\underline{\lam}\omega,\   \ &0<x<L,\\[2mm]
d_I\omega_x=\F{q}{2}\omega,\ &x=0,L.
\end{array}
\right.
\eeq
It is obvious that $\lam_1(d_I,q,\gam-\s\be)\ge\underline{\lam}$.
Thus,
we will only  prove that $\underline{\lam}$ is bounded from below.
The holding for \emph{Step 1} then follows directly.

Notice that $\underline{\lam}$ admits the following
variational characterization:
\beq
\underline{\lam}=
\inf\limits_{\omega\in H^1((0,L)),\omega\neq0}
\F{d_I\i_0^L|\omega_x-\f{q}{2d_I}\omega|^2\dx+\i_0^L(\gam_*-\be^*\s(x))\omega^2\dx}{\i_0^L\omega^2\dx}.
\eeq

Take $w_{d_S}$ the corresponding positive principal eigenfunction,
normalized by $\|w_{d_S}\|_2=1$,
for all $d_S>0$.
Suppose on the contrary that there exists a subsequence of $d_S\rightarrow0$ (we still denote it by itself),
such that
\beq
\label{r1->infinity}
\lim\limits_{d_S\rightarrow0}\underline{\lam}=-\infty.
\eeq
Now we have some claims.

\emph{Claim 1.1.}
$\s_x>0$ on $[0,L]$ for all $d_S>0$,
provided that
$\mu'(x)\le0,\Lam'(x)\ge0$.

Let $a(x)=\F{\s_x}{\s}$.
Hence it follows from \eqref{hattilde S} that
\beq
\left\{\arraycolsep=1.5pt
\begin{array}{ll}
d_S a_{xx}+(2d_S a-q)a_x-\Lam(x)a=\mu'(x)-\F{\Lam'(x)}{\s},\   \ &0<x<L,\\[2mm]
a(0)=a(L)=\F{q}{d_S}>0,\ &x=0,L,
\end{array}
\right.
\eeq
which together with the maximum principle yields that $a>0$,
and hence \emph{Claim 1.1} is correct.

\emph{Claim 1.2.}
$(w_{d_S})_x>0$ on $[0,L]$ for all $d_S>0$.

Since $(w_{d_S})_x(0)>0,(w_{d_S})_x(L)>0$,
then there exist $0<x_1<x_2<L$ such that $x_1$($x_2$) is a local maximum(minimum) point of $w_{d_S}$ in $(0,L)$,
if $w_{d_S}$ is not monotonically increasing on $[0,L]$ for some $d_S>0$.
Furthermore, we reach that
$(w_{d_S})_{xx}(x_1)\le0$ and $(w_{d_S})_{xx}(x_2)\ge0$.
By the first equation of \eqref{new eigenvalue problem},
we have that
\beq
\F{q^2}{4d_I}+\gam_*-\be^*\s(x_1)-\underline{\lam}\le0\le\F{q^2}{4d_I}+\gam_*-\be^*\s(x_2)-\underline{\lam}.
\eeq
This is impossible due to $x_1<x_2$ and \emph{Claim 1.1}.
Thus, $w_{d_S}$ is increasing on $[0,L]$ for all $d_S>0$.

\emph{Claim 1.3.}
$w_{d_S}\rightarrow0$ locally uniformly in $[0,L)$ as $d_S\rightarrow0$.

Owing to \emph{Claim 1.2},
we only to prove for any $x\in[0,L)$ there holds $\lim\limits_{d_S\rightarrow0}w_{d_S}(x)=0$.
Without loss of generality,
due to \emph{Claim 1.2},
suppose now that there exists some $x_0\in(0,L)$ such that $\lim\limits_{d_S\rightarrow0}w_{d_S}(x_0)>0$.
Choosing a nonnegative text function $v\in C_0^2([0,L])$ with
\beq
\begin{cases}
v>0\ \text{in}\ (\F{x_0}{2},\F{x_0+L}{2}),\\
v\equiv0\ \text{on}\ [0,\F{x_0}{2}]\cup[\F{x_0+L}{2},L],
\end{cases}
\eeq
and multiplying both sides of \eqref{new eigenvalue problem} by $v$,
it indicates from the integrating by parts that
\beq
\bea{l}
\underline{\lam}=
\Big|\F{-d_I\io w_{d_S}v_{xx}\dx+\io(\f{q^2}{4d_I}+\gam_*-\be^*\s) w_{d_S}v\dx}
{\io w_{d_S}v\dx}\Big|\\
\ge\F{|-d_I\i_{\f{x_0}{2}}^{\f{x_0+L}{2}} w_{d_S}v_{xx}\dx+\i_{\f{x_0}{2}}^{\f{x_0+L}{2}}(\f{q^2}{4d_I}+\gam_*-\be^*\s) w_{d_S}v\dx|}
{\i_{\f{x_0}{2}}^{\f{x_0+L}{2}}w_{d_S}v\dx}.
\eea
\eeq
From the choice of $v$,
we can get from the inequality above that $\lim\limits_{d_S\rightarrow0}\underline{\lam}>0$,
which contradicts to \eqref{r1->infinity}.
Hence, \emph{Claim 1.3} is true.

Given $\emph{l}\in(0,L)$,\emph{Claim 1.3} yields that $\lim\limits_{d_S\rightarrow0}w_{d_S}(L-\emph{l})=0$.
Notice that $\|w_{d_S}\|_2=1$,
$\lim\limits_{d_S\rightarrow0}w_{d_S}(L)=+\infty$ can be deduced by \emph{Claim 1.1} and \emph{Claim 1.2}.
Hence,
\beq
\label{2}
\bea{l}
\Io|(w_{d_S})_x|^2\dx\ge\I_{L-\emph{l}}^{L}|(w_{d_S})_x|^2\dx
\ge\F{(\i_{L-\emph{l}}^{L}(w_{d_S})_x\dx)^2}{\emph{l}}
=\F{(w_{d_S}(L)-w_{d_S}(L-\emph{l}))^2}{\emph{l}}\\
=\F{(w_{d_S}(L))^2-2w_{d_S}(L)w_{d_S}(L-\emph{l})+(w_{d_S}(L-\emph{l}))^2}{\emph{l}}\\
\ge\F{2(w_{d_S}(L))^2+2(w_{d_S}(L-\emph{l}))^2}{\emph{l}}\rightarrow+\infty,\ \text{as}\ d_S\rightarrow0,
\eea
\eeq
where we have used the H\"{o}lder inequality and the Cauchy's inequality.
Thus, without loss of generality,
we may assume that
\beq
\Io|(w_{d_S})_x|^2\dx\ge\F{4q^2}{d_I^2},
\eeq
for all small $d_S$.
Therefore, by \eqref{2},
we obtain that
\beq
\label{3}
\bea{l}
\Io((w_{d_S})_x-\F{q}{2d_I}w_{d_S})^2\dx
=\Io|(w_{d_S})_x|^2\dx-\F{q}{d_I}\Io w_{d_S}(w_{d_S})_x\dx+\F{q^2}{4d_I^2}\\
\ge\Io|(w_{d_S})_x|^2\dx-\F{q}{d_I}(\Io |(w_{d_S})_x|^2\dx)^{\f12}+\F{q^2}{4d_I^2}\\
\ge\F12\Io|(w_{d_S})_x|^2\dx+\F{q^2}{4d_I^2}
\ge\F{(w_{d_S}(L))^2+(w_{d_S}(L-\emph{l}))^2}{\emph{l}}+\F{q^2}{4d_I^2}.
\eea
\eeq

Define $T(x)=\I_0^{x}\s(y)\dy$.
Then, integrating by parts over $(0,L)$,
together with \emph{Claim 1.2},
yields that
\beq
\Io\s(x)w_{d_S}^2(x)\dx\le w_{d_S}^2(L)T(L)-2\Io T(x)w_{d_S}(x)(w_{d_S})_x(x)\dx
\le \F{w_{d_S}^2(L)\|\Lam\|_1}{\mu_*}.
\eeq
Therefore,
by applying the inequality above,
\eqref{3} and the definition of $w_{d_S}$,
it follows that
\beq
\bea{l}
\underline{\lam}=
d_I\I_0^L((w_{d_S})_x-\f{q}{2d_I}w_{d_S})^2\dx
+\I_0^L(\gam_*-\be^*\s)(w_{d_S})^2\dx\\
\ge\F{d_I(w_{d_S}(L))^2+d_I(w_{d_S}(L-\emph{l}))^2}{\emph{l}}+\F{q^2}{4d_I}+\gam_*-\F{w^2_{d_S}(L)\be^*\|\Lam\|_1}{\mu_*}\\
\rightarrow+\infty\ \text{as}\ d_S\rightarrow0,
\eea
\eeq
provided that \emph{l} is small enough.
This contradicts to \eqref{r1->infinity}.
Hence we have verified \emph{Step 1}.

\emph{Step 2.}
For any given $\eps\in(0,L)$,
there holds
$\lim\limits_{d_S\rightarrow0}\s(L-\eps)=\s^{\infty}(L-\eps)$,
where $\s^{\infty}$ satisfies \eqref{s infty}.

Note that \eqref{9} implies that $\s$ is locally uniformly bounded in $[0,L)$ for $d_S\rightarrow0$.
For $\eps>0$ small enough,
after passing to a subsequence of $d_S\rightarrow0$,
we may assume that
$\s\rightarrow\s^{\infty}\ge0$ weakly in $L^2((0,L-\F{\eps}{2}))$
for some $\s^{\infty}\in L^2((0,L-\F{\eps}{2}))$.
Since
$\Io\s(x)\dx\le\F{\|\Lam\|_1}{\mu_*}$, it then indicates from a diagonal argument(up to a further subsequence of $d_S\rightarrow0$) that
\beq
\Io\s^{\infty}(x)\dx
=\lim\limits_{\eps\rightarrow0}
\I_0^{L-\f{\eps}{2}}\s^{\infty}(x)\dx
=\lim\limits_{\eps\rightarrow0}
\lim\limits_{d_S\rightarrow0}
\I_0^{L-\f{\eps}{2}}\s(x)\dx\le\F{\|\Lam\|_1}{\mu_*}.
\eeq
Hence $\s^{\infty}\in L^1((0,L))$.

Take any test function $\xi\in C_0^2([0,L-\F{\eps}{2}])$.
Multiplying the both sides of \eqref{hattilde S} by $\xi$ and integrating by parts in $(0,L-\F{\eps}{2})$,
there holds \beq
\begin{array}{l}
d_S\I_0^{\f{\eps}{2}}\s\xi_{xx}\dx+q\I_0^{\f{\eps}{2}}\s\xi_x\dx+\I_0^{\f{\eps}{2}}(\Lam(x)-\mu(x)\s(x))\xi\dx=0.
\end{array}
\eeq
Hence letting $d_S\rightarrow0$ and due to the arbitrariness of $\eps$,
we have that $\s^{\infty}$ is a weak and
then a classical solution of the first equation in \eqref{s infty}.

Set $\s(x_{d_S})=\min\limits_{x\in[0,L]}\s$.
Owing to $\s_x(L)>0$,
it is clear that $x_{d_S}\in[0,L)$.

\emph{Claim 2.1}
$\lim\limits_{d_S\rightarrow0}\s(x_{d_S})=0$.

If not, then there exists some $\del>0$ and $0<d_{S_n}<\F{1}{n}$
such that $\s(x_{d_{S_n}})\ge\del$ on $[0,L]$.
Let $w(x)=e^{-\f{q}{d_{S_n}}(1-C\f{d_{S_n}}{q^2})x}\s(x)$.
Therefore $w$ satisfies
\beq
\label{6}
\left\{\arraycolsep=1.5pt
\begin{array}{ll}
d_{S_n} w_{xx}+q(1-2C\F{d_{S_n}}{q^2})w_x+\big[-C(1-C\F{d_{S_n}}{q^2})-\mu(x)\big]w+\Lam e^{-\f{q}{d_{S_n}}(1-C\f{d_{S_n}}{q^2})x}=0,\   \ &0<x<L,\\[2mm]
w_x=\F{C}{q}w>0,\ &x=0,L.
\end{array}
\right.
\eeq
Take $w(\hat{x})=\max\limits_{x\in[0,L]}w(x)$.
Clearly, $\hat{x}>0$,
since $w_x(0)>0$.
If $\hat{x}\in(0,L)$,
then $w_x(\hat{x})=0$ and $w_{xx}(\hat{x})\le0$.
By applying \eqref{6},
one can have that
\beq
\big(-C(1-C\F{d_{S_n}}{q^2})-\mu(\hat{x})\big)w(\hat{x})+\F{\Lam(\hat{x})}{\s(\hat{x})}w(\hat{x})\ge0.
\eeq
However,
\beq
\bea{l}
\big(-C(1-C\F{d_{S_n}}{q^2})-\mu(\hat{x})\big)w(\hat{x})+\F{\Lam(\hat{x})}{\s(\hat{x})}w(\hat{x})\\
\le
\big(-C(1-C\F{d_{S_n}}{q^2})+\F{\Lam(\hat{x})}{\del}\big)w(\hat{x})<0,
\eea
\eeq
provided that $C>\F{\Lam^*}{\del}$
and $n$ is large enough.
Then a contradiction occurs,
and therefore $\hat{x}=L$.
In turn, $w(x)\le w(L)$ for all $x\in[0,L]$.
So, $\s(x)\le\s(L)e^{-\f{q}{d_{S_n}}(1-C\f{d_{S_n}}{q^2})(L-x)}$.
Due to \eqref{s(L)<=},
one can see that
this is a contradiction to our assumption that $\s(x_{d_{S_n}})\ge\del>0$
and thus we have
proved \emph{Claim 2.1}.

If $x_{d_S}\in(0,L)$,then $\s_{xx}(x_{d_S})\ge0$
and $\s_x(x_{d_S})=0$.\eqref{hattilde S} hence yields that
\beq
\Lam_*-\mu(x_{d_S})\s(x_{d_S})\le
\Lam(x_{d_S})-\mu(x_{d_S})\s(x_{d_S})\le0,
\eeq
which contradicts to \emph{Claim 2.1}.
Thus, $x_{d_S}=0$ for all small $d_S>0$.
As a
consequence,
we have shown $\s^{\infty}(0)=0$.
Therefore, $\s^{\infty}$ is the classical solution of \eqref{s infty}
and $\lim\limits_{d_S\rightarrow0}\s(L-\eps)=\s^{\infty}(L-\eps)$.

\emph{Step 3.}
For any given $\eps\in(0,L)$,
there holds
$\lim\limits_{d_S\rightarrow0}d_S\s_x(L-\eps)=0$.

Multiplying the first equation of \eqref{hattilde S} by $e^{-\f{qx}{d_S}}$ and integrating over $(x,L)$,
we get that
\beq
d_S\s_x e^{-\f{qx}{d_S}}-d_S\s_x(L)e^{-\f{qL}{d_S}}\le\F{d_S\Lam^*}{q}(e^{-\f{qx}{d_S}}-e^{-\f{qL}{d_S}}),
\eeq
which implies that
\beq
\label{7}
\s_x(x)\le\s_x(L)e^{-\f{q(L-x)}{d_S}}+\F{\Lam^*}{q}(1-e^{-\f{q(L-x)}{d_S}})\ \text{for any}\ x\in[0,L].
\eeq
Again, multiplying the first equation of \eqref{hattilde S} by $e^{-\f{qx}{d_S}}$ and integrating over $(x,L)$ yield that
\beq
d_S\s_x(L) e^{-\f{qL}{d_S}}-d_S\s_x(x)e^{-\f{qx}{d_S}}\le\mu^*\I_x^L\s(y)e^{-\f{qy}{d_S}}\dy.
\eeq
Then
\beq
\s_x(L) e^{-\f{q(L-x)}{d_S}}\le \s_x(x)+\F{\mu^*}{d_S}e^{\f{qx}{d_S}}\I_x^L\s(y)e^{-\f{qy}{d_S}}\dy.
\eeq
We now integrate the above inequality from $0$ to $L$
to arrive at
\beq
\bea{l}
\s_x(L)\F{d_S}{q}(1-e^{-\f{qL}{d_S}})\le\s(L)
+\F{\mu^*}{d_S}\Io(e^{\f{qx}{d_S}}\I_x^L\s(y)e^{-\f{qy}{d_S}}\dy)\dx\\
\le\s(L)+\F{\mu^*}{q}\Io\s(x)\dx
\le\s(L)+\F{\|\Lam\|_1\mu^*}{q\mu_*},
\eea
\eeq
where we have used integrating by parts.
Then we get that
\beq
\s_x(L)\le\F{q}{d_S(1-e^{-\f{qL}{d_S}})}(\s(L)+\F{\|\Lam\|_1\mu^*}{q\mu_*}).
\eeq
Combined with \eqref{s(L)<=} and \eqref{7},
it then follows that $\s_x$ is locally uniformly bounded in $[0,L)$ for $d_S\rightarrow0$.
In particular,
$\s_x(L-\eps)<+\infty$ for $d_S\rightarrow0$.
Therefore the proof of \emph{Step 3} is finished.

\emph{Step 4.}
$\lim\limits_{d_S\rightarrow0}\lam_1(d_I,q,\gam-\s\be)=\bar{\lam}$.

Let
$\psi(x)>0$ be the principal eigenfunction of \eqref{lam1(d q ga-s)}, normalized by $\Io e^{-\f{qx}{d_I}}\psi(x)\dx=1$.
Upon the transformation $\phi=e^{-\f{qx}{d_I}}\psi$,
we have that
\beq
\left\{\arraycolsep=1.5pt
\begin{array}{ll}
-d_I\phi_{xx}-q\phi_x+(\gam(x)-\be(x)\s(x))\phi=\lam_1\phi,\   \ &0<x<L,\\[2mm]
\phi_x=0,\ &x=0,L,\\
\Io\phi\dx=1.&
\end{array}
\right.
\eeq
Since $\Io\be(x)\s(x)\dx\le\be^*\Io\s(x)\dx\le\F{\be^*\|\Lam\|_1}{\mu_*}$,
and since \eqref{lam1 variational} leads to
\beq
\lam_1\le\F{\io(\gam(x)-\be(x)\s^{\infty}(x))e^{\f{qx}{d_I}}\dx}{\io e^{\f{qx}{d_I}}\dx}
\le\gam^*,
\eeq
it then follows from the Harnack inequality that
\beq
\max\limits_{0\le x\le L}\phi\le C\min\limits_{0\le x\le L}\phi.
\eeq
Here $C$ is independent of $d_S>0$.
Moreover,
thanks to
\beq
\min\limits_{0\le x\le L}\phi\le\F{1}{L}\Io\phi\dx\le\F{1}{\sqrt{L}}(\Io\phi^2\dx)^{\f12}\le\F{1}{\sqrt{L}},
\eeq
we have that
$\phi$ is uniformly bounded
and furthermore $\psi$ is also uniformly bounded in $[0,L]$.
We now integrate \eqref{lam1(d q ga-s)} over $[0,x]$ to obtain
\beq
d_I\psi_x=q\psi+\I_0^x(\gam(y)-\be(y)\s(y)-\lam_1)\psi(y)\dy
\le q\psi+\I_0^L(\gam(y)-\lam_1)\psi(y)\dy.
\eeq
Then one observes that $\psi_x$ is also uniformly
bounded over $[0,L]$ and furthermore $\psi$ is bounded in $C^1([0,L])$ independent of $d_S>0$.
Therefore it can be checked from \eqref{lam1(d q ga-s)} that $\psi_{xx}$ is locally uniformly
bounded in $[0,L)$.
Now, passing to a subsequence of $d_S\rightarrow0$,
denoted by itself for simplicity,
we can assume
that $\lim\limits_{d_S\rightarrow0}\lam_1=\bar{\lam}$
and $\psi\rightarrow\bar{\psi}\ge,\not\equiv0$ in $C([0,L])\cap C^1_{loc}([0,L))$ as $d_S\rightarrow0$.
Indeed, $(\bar{\lam},\bar{\psi})$
satisfies the first equation of \eqref{lim eigen problem} (in the weak sense) and
$d_I\bar{\psi}_x(0)-q\bar{\psi}(0)=0$.

Finally, we determine the boundary condition of $\bar{\psi}$ at $x=L$.
We fix any small $\eps\in(0,L)$ and
integrate \eqref{lam1(d q ga-s)} over $[L-\eps,L]$ to deduce that
\beq
d_I\psi_x(L-\eps)-q\psi(L-\eps)+\I_{L-\eps}^L(\gam(x)-\be(x)\s(x))\psi(x)\dx
=\lam_1\I_{L-\eps}^L\psi(x)\dx.
\eeq
Combined with the mean value theorem,
we arrive at
\beq
d_I\psi_x(L-\eps)-q\psi(L-\eps)+\I_{L-\eps}^L\gam(x)\psi(x)\dx
-\be(x_{\eps})\psi(x_{\eps})\I_{L-\eps}^L\s(x)\dx=\lam_1\I_{L-\eps}^L\psi(x)\dx,
\eeq
for some $x_{\eps}\in(L-\eps,L)$.
Notice that integrating \eqref{hattilde S} over $[L-\eps,L]$
and the mean value theorem
give rise to
\beq
\I_{L-\eps}^L\s(x)\dx
=\F{-d_S\s_x(L-\eps)+q\s(L-\eps)+\I_{L-\eps}^L\Lam(x)\dx}{\mu(x^\eps)},
\eeq
for some $x^\eps\in(L-\eps,L)$.
Because of \emph{Step 2} and \emph{Step 3},
one can obtain that
\beq
d_I\bar{\psi}_x(L)-q\bar{\psi}(L)=\F{\be(L)\bar{\psi}(L)q\s^{\infty}(L)}{\mu(L)},
\eeq
by letting $d_S\rightarrow0$ first and then $\eps\rightarrow0$.
Since
integrating \eqref{s infty} from $0$ to $L$ leads to $q\s^{\infty}(L)=\|\Lam\|_1-\Io\mu(x)\s^{\infty}(x)\dx$,
it hence deduced that
\beq
d_I\bar{\psi}_x(L)-q\bar{\psi}(L)=\be(L)\bar{\psi}(L)\mu(L)N_S,
\eeq
where $N_S$ is as in \eqref{ns}.
The proof of
Theorem \ref{lim of lam1} is then finished.
\proofend

The following lemma is about threshold-type dynamics in terms of the
basic reproduction number.
Since Lemma \ref{stable unstable} can be proved by the similar method of Lemma 2.4 in \cite{c},
and of Theorem 3.3 in \cite{sc} for a similar
problem,
so the details are omitted here.

\begin{lem}
\label{stable unstable}
If $\rr_0<1$,
then the DFE is linearly stable,
but if $\rr_0>1$ then it is linearly unstable
and there exists a constant $\eta>0$ independent of the initial data,
such that any solution $(S,I)$ of \eqref{P} satisfies
\beq
\liminf\limits_{t\rightarrow+\infty}S(x,t)\ge\eta\ \text{and}\  \liminf\limits_{t\rightarrow+\infty}I(x,t)\ge\eta\ \text{uniformly for}\ x\in[0,L],
\eeq
and hence, the disease persist uniformly. Furthermore, \eqref{P} admits at least one EE when $\rr_0>1$.
\end{lem}

\section{Global attractivity of the EE}
\label{Global attractivity of the EE}

In this section,
we shall study the global attractivity of the EE of \eqref{P} by assuming that
$\Lam(x)\be(x)>\gam(x)\mu(x)$,
and both $\F{\Lam(x)}{\mu(x)e^{\f{qx}{d_S}}}$
and $\F{\Lam(x)\be(x)-\gam(x)\mu(x)}{\gam(x)\mu(x)e^{\f{qx}{d_I}}}$ are positive constants.
Upon the transform $u=Se^{-\f{qx}{d_S}}$, $v=Ie^{-\f{qx}{d_I}}$,
it deduces that the system \eqref{P}
is equivalent
to the following one
\beq
\label{uv}
\left\{\arraycolsep=1.5pt
\begin{array}{lll}
u_t=d_S u_{xx}+q u_x+\Lam(x)e^{-\f{qx}{d_S}}-\mu(x)u-\beta(x)\F{uve^{\f{qx}{d_I}}}{1+mve^{\f{qx}{d_I}}}+\gam(x)ve^{(\f{q}{d_I}-\f{q}{d_S})x},\   \ &0<x<L,\ t>0,\\[2mm]
v_t=d_I v_{xx}+q v_x+\beta(x)\F{uve^{\f{qx}{d_S}}}{1+mve^{\f{qx}{d_I}}}-\gam(x)v, \   \ &0<x<L,\ t>0,\\[2mm]
u_x=v_x=0,\ &x=0,L,\ t>0,\\[2mm]
u(x,0)=S_0(x)e^{-\f{qx}{d_S}}\ge,\not\equiv0,\ v(x,0)=I_0(x)e^{-\f{qx}{d_I}}\ge,\not\equiv0, \ &0<x<L.
\end{array}
\right.
\eeq
Hence $(\F{\Lam(x)}{\mu(x)e^{\f{qx}{d_S}}},\F{\Lam(x)\be(x)-\gam(x)\mu(x)}{\gam(x)\mu(x)e^{\f{qx}{d_I}}})$ is a positive constant EE of \eqref{uv}.

\two
Denote $(\F{\Lam(x)}{\mu(x)e^{\f{qx}{d_S}}},\F{\Lam(x)\be(x)-\gam(x)\mu(x)}{\gam(x)\mu(x)e^{\f{qx}{d_I}}})$ by $(\bar{u},\bar{v})$ for convenience.
We construct the following Lyapunov functional
\beq
F(t)=\F12\Io e^{\f{qx}{d_S}}(u-\bar{u})^2\dx+\F12\Io e^{\f{qx}{d_I}}(v-\bar{v})^2\dx.
\eeq
By some calculations and integrating by parts,
we get
that
\beq
\label{1}
\arraycolsep=1.5pt
\begin{array}{l}
F'(t)=\Io(u-\bar{u})\big[(e^{\f{qx}{d_S}}u_x)_x+\Lam(x)-\mu(x)ue^{\f{qx}{d_S}}-\be(x)e^{(\f{q}{d_I}+\f{q}{d_S})x}\F{uv}{1+mve^{\f{qx}{d_I}}}+\gam(x)ve^{\f{qx}{d_I}}\big]\dx\\
+\Io(v-\bar{v})\big[(e^{\f{qx}{d_I}}v_x)_x+\be(x)e^{(\f{q}{d_I}+\f{q}{d_S})x}\F{uv}{1+mve^{\f{qx}{d_I}}}-\gam(x)ve^{\f{qx}{d_I}}\big]\dx\\
=-\Io e^{\f{qx}{d_S}}u_x^2\dx-\Io e^{\f{qx}{d_I}}v_x^2\dx\\
+\Io(u-\bar{u})\big[\Lam(x)-\mu(x)ue^{\f{qx}{d_S}}-\be(x)e^{(\f{q}{d_I}+\f{q}{d_S})x}\F{uv}{1+mve^{\f{qx}{d_I}}}+\gam(x)ve^{\f{qx}{d_I}}\big]\dx\\
+\Io(v-\bar{v})\big[\be(x)e^{(\f{q}{d_I}+\f{q}{d_S})x}\F{uv}{1+mve^{\f{qx}{d_I}}}-\gam(x)ve^{\f{qx}{d_I}}\big]\dx.
\end{array}
\eeq
Since the definition of $(\bar{u},\bar{v})$ yields that
\beq
\left\{\arraycolsep=1.5pt
\begin{array}{l}
\Lam(x)e^{-\f{qx}{d_S}}-\mu(x)\bar{u}=0,\\[2mm]
\beta(x)\F{\bar{u}e^{\f{qx}{d_S}}}{1+m\bar{v}e^{\f{qx}{d_I}}}-\gam(x)=0,
\end{array}
\right.
\eeq
we then obtain that
\beq
\bea{l}
F'(t)=-\Io e^{\f{qx}{d_S}}u_x^2\dx-\Io e^{\f{qx}{d_I}}v_x^2\dx\\
+\Io e^{\f{qx}{d_S}}\mu(x)(u-\bar{u})^2\dx
-\Io\be(x)e^{(\f{q}{d_S}+\f{q}{d_I})x}\F{v(u-\bar{u})^2}{1+mve^{\f{qx}{d_I}}}\dx\\
+\Io\be(x)e^{(\f{q}{d_S}+\f{q}{d_I})x}\F{vm\bar{u}e^{\f{qx}{d_I}}}{(1+mve^{\f{qx}{d_I}})(1+m\bar{v}e^{\f{qx}{d_I}})}(u-\bar{u})(v-\bar{v})\dx\\
+\Io\be(x)e^{(\f{q}{d_S}+\f{q}{d_I})x}\F{v(u-\bar{u})(v-\bar{v})}{1+mve^{\f{qx}{d_I}}}\dx
-\Io\be(x)e^{(\f{q}{d_S}+\f{q}{d_I})x}\F{\bar{u}(v-\bar{v})^2me^{\f{qx}{d_I}}}{(1+mve^{\f{qx}{d_I}})(1+m\bar{v}e^{\f{qx}{d_I}})}\dx\\
:=-\Io e^{\f{qx}{d_S}}u_x^2\dx-\Io e^{\f{qx}{d_I}}v_x^2\dx
+\Io(A(u-\bar{u})^2+B(u-\bar{u})(v-\bar{v})+C(v-\bar{v})^2)\dx,
\eea
\eeq
where $A<0$ and $C<0$ can be derived by some direct calculations.
Since
\beq
\lim\limits_{m\rightarrow+\infty}\F{(B^2-4AC)(1+mve^{\f{qx}{d_I}})(1+m\bar{v}e^{\f{qx}{d_I}})}{\be ve^{2(\f{q}{d_S}+\f{q}{d_I})x}}=-\infty,
\eeq
there exists a $M>0$ such that for any $m>M$,
$B^2-4AC<0$.
Consequently,
\beq
F'(t)\le-\Io e^{\f{qx}{d_S}}u_x^2\dx-\Io e^{\f{qx}{d_I}}v_x^2\dx-\Io\del v\big((u-\bar{u})^2+(v-\bar{v})^2\big)\dx\le0
\eeq
for some $0<\del<<1$ and for any $m>M$.
Then we get that $F'(t)=0\Rightarrow(u,v)=(\bar{u},\bar{v})$.
By the LaSalle invariance principle,
there holds $\lim\limits_{t\rightarrow+\infty}(u(x,t),v(x,t))=(\bar{u},\bar{v})$.
This then implies that the EE $(\bar{S},\bar{I})$ is globally asymptotically
stable for \eqref{P}.
\proofend

\section{Asymptotic profiles of the EE}
\label{Asymptotic profiles of the EE}

It follows from Lemma \ref{stable unstable} that \eqref{P} admits at least one EE when $\rr_0>1$.
In this section,
we consider the asymptotic profiles of the EE in the following six cases: (i): $q\rightarrow+\infty$;
(ii)
$d_S\rightarrow0$;
(iii)
$d_I\rightarrow0$;
(iv) $d_S\rightarrow+\infty$;
(v)$d_I\rightarrow+\infty$;
(vi)$m\rightarrow+\infty$.
Clearly, the EE satisfies \eqref{equilibrium}.
%the following system:
%\beq
%\label{EE's equation}
%\left\{\arraycolsep=1.5pt
%\begin{array}{ll}
%d_S S_{xx}-q S_x+\Lam(x)-\mu(x)S-\beta(x)\F{SI}{1+mI}+\gam(x)I=0,\   \ &0<x<L,\\[2mm]
%d_I I_{xx}-q I_x+\beta(x)\F{SI}{1+mI}-\gam(x)I=0, \   \ &0<x<L,\\[2mm]
%d_S S_x-q S=d_I I_x-q I=0,\ &x=0,L.
%\end{array}
%\right.
%\eeq

\subsection{A priori estimates of solution of \eqref{P}}

We begin with the prior estimates
for solutions to \eqref{P},
which will be used in obtaining the limiting profiles of the EE.

\begin{lem}
For any $d_S,d_I,q>0$,
any positive solution $(S,I)$ of \eqref{equilibrium} satisfies the following estimates:

(1)
\beq
\label{estimate1}
\F{\|\Lam\|_1}{\mu^*}\le\Io S(x)\dx\le \F{\|\Lam\|_1}{\mu_*},\ \Io I(x)\dx\le\F{\be^*\|\Lam\|_1}{m\mu_*\gam_*};
\eeq

(2)
\beq
\label{s<=}
S(x)\le\F{\|\Lam\|_1}{q}(1+\F{\be^*}{m\mu_*})(1-e^{-\f{q}{d_S}(L-x)})
+\F{\|\Lam\|_1e^{-\f{q(L-x)}{d_S}}}{d_S(1-e^{-\f{qL}{d_S}})}\big[\F{q}{\mu_*}+(1+\F{\be^*}{m\mu_*})(L-\F{d_S}{q}(1-e^{-\f{qL}{d_S}}))\big],
\eeq
for all $x\in[0,L]$;

(3)
\beq
\label{i<=}
I(x)\le\F{\|\Lam\|_1\be^*}{qm\mu_*}(1-e^{-\f{q}{d_I}(L-x)})
+\F{\be^*\|\Lam\|_1e^{-\f{q(L-x)}{d_I}}}{d_I m\mu_*(1-e^{-\f{qL}{d_I}})}\big[\F{q}{\gam_*}+(L-\F{d_I}{q}(1-e^{-\f{qL}{d_I}}))\big],
\eeq
for all $x\in[0,L]$.
\end{lem}

\proof
(1) We integrate both sides of the first two equations of \eqref{equilibrium} in $(0,L)$ to have that
\beq
\label{10}
\begin{cases}
\Io\Lam(x)\dx=\Io\mu(x)S(x)\dx,\\
\Io\gam(x)I(x)\dx=\Io\F{\be(x)SI}{1+mI}\dx<\F{1}{m}\Io\be(x)S(x)\dx.
\end{cases}
\eeq
Hence, the desired estimates \eqref{estimate1} can be deduced immediately.

(2) Integrating the equation of $S(x)$ in \eqref{equilibrium} from $x$ to $L$ yields
\beq
\bea{ll}
-d_S S_x+qS&=\I_x^L(-\Lam(y)+\mu(y)S+\F{\be(y)SI}{1+mI})\dy\\
&<\I_x^L(\mu(y)S+\F{\be(y)SI}{1+mI})\dy\le\Io\mu(y)S\dy+\Io\F{\be(y)SI}{1+mI})\dy\\
&=\Io\Lam(y)\dy+\Io\gam(y)I(y)\dy\le\|\Lam\|_1(1+\F{\be^*}{m\mu_*}).
\eea
\eeq
By the similar arguments as in \eqref{8}-\eqref{9},
we conclude that
\beq
\label{a}
\begin{cases}
S(L)\le\F{1}{d_S(1-e^{-\f{qL}{d_S}})}\big[\F{\|q\Lam\|_1}{\mu_*}+\|\Lam\|_1(1+\F{\be^*}{m\mu_*})(L-\F{d_S}{q}(1-e^{-\f{qL}{d_S}}))\big],\\
S(x)\le\F{\|\Lam\|_1}{q}(1+\F{\be^*}{m\mu_*})(1-e^{-\f{q}{d_S}(L-x)})
+\F{\|\Lam\|_1e^{-\f{q(L-x)}{d_S}}}{d_S(1-e^{-\f{qL}{d_S}})}\big[\F{q}{\mu_*}+(1+\F{\be^*}{m\mu_*})(L-\F{d_S}{q}(1-e^{-\f{qL}{d_S}}))\big].
\end{cases}
\eeq
Actually,
some similar calculations can imply that
\beq
\label{ix<= iL<=}
\begin{cases}
I(L)\le\F{\be^*\|\Lam\|_1}{d_Im\mu_*(1-e^{-\f{qL}{d_I}})}\big[\F{q}{\gam_*}+(L-\F{d_I}{q}(1-e^{-\f{qL}{d_I}}))\big],\\
I(x)\le
\F{\|\Lam\|_1\be^*}{qm\mu_*}(1-e^{-\f{q}{d_I}(L-x)})
+\F{\be^*\|\Lam\|_1e^{-\f{q(L-x)}{d_I}}}{d_I m\mu_*(1-e^{-\f{qL}{d_I}})}\big[\F{q}{\gam_*}+(L-\F{d_I}{q}(1-e^{-\f{qL}{d_I}}))\big].
\end{cases}
\eeq
\proofend

\subsection{Concentration behaviors of EE as $q\rightarrow+\infty$}

This subsection is concerned with the asymptotic profile of EE,
if exists,
with respect to large advection
rate.

\three
We split the proof into three steps.

\emph{Step 1.}
Any positive solution $(S(x),I(x))$ of \eqref{equilibrium} satisfies
\beq
(S(x),I(x))\rightarrow(0,0)\ \text{locally uniformly in}\ [0,L)\ \text{as}\ q\rightarrow+\infty.
\eeq
In fact, this statement follows from estimates \eqref{s<=} and \eqref{i<=}.

\emph{Step 2.}
Convergence of $(S(L),I(L))$

Upon the transform $(a(y),b(y)):=(\F{1}{q}S(L-\F{y}{q}),\F{1}{q}I(L-\F{y}{q}))$,
there holds
\beq
\label{a b's equ}
\left\{\arraycolsep=1.5pt
\begin{array}{lll}
d_S a_{yy}+a_y+\F{1}{q^3}\Lam(L-\F{y}{q})-\F{1}{q^2}\mu(L-\F{y}{q})a-\F{ab}{q(1+mqb)}\be(L-\F{y}{q})+\F{b}{q^2}\gam(L-\F{y}{q})=0,\   \ &0<y<qL,\\[2mm]
d_I b_{yy}+b_y+\F{ab}{q(1+mqb)}\be(L-\F{y}{q})-\F{b}{q^2}\gam(L-\F{y}{q})=0, \   \ &0<y<qL,\\[2mm]
d_S a_y+a=d_I b_y+b=0,\ &y=0,qL.
\end{array}
\right.
\eeq
Using the estimates \eqref{s<=} and \eqref{i<=} yield that
\beq
\label{a<= b<=}
\begin{cases}
a(y)\le\F{\|\Lam\|_1}{q^2}(1+\F{\be^*}{m\mu_*})(1-e^{-\f{y}{d_S}})
+\F{\|\Lam\|_1e^{-\f{y}{d_S}}}{d_S(1-e^{-\f{qL}{d_S}})q}\big[\F{q}{\mu_*}+(1+\F{\be^*}{m\mu_*})(L-\F{d_S}{q}(1-e^{-\f{qL}{d_S}}))\big],\\
b(y)\le\F{\|\Lam\|_1\be^*}{q^2m\mu_*}(1-e^{-\f{y}{d_I}})
+\F{\be^*\|\Lam\|_1e^{-\f{y}{d_I}}}{d_I m\mu_*q(1-e^{-\f{qL}{d_I}})}\big[\F{q}{\gam_*}+(L-\F{d_I}{q}(1-e^{-\f{qL}{d_I}}))\big],
\end{cases}
\eeq
for any $0\le y\le qL$.
Then for any fixed positive constant $M$ and all $q>\F{M}{L}$,
$a(y)$ and
$b(y)$ are uniformly bounded on $[0,M]$.
Consequently,
the standard $L^p$-estimates of elliptic equations and Sobolev embedding theorem yield that
$\|a\|_{C^{1,\al}([0,M])}\le C$,
$\|b\|_{C^{1,\al}([0,M])}\le C$ for some $\al\in(0,1)$.
Passing to a subsequence if necessary,
$\{q_n\}$ satisfies $q_n\rightarrow+\infty$ as $n\rightarrow+\infty$,
then there further hold
\beq
\lim\limits_{n\rightarrow+\infty}(a_n(y),b_n(y)):=\lim\limits_{n\rightarrow+\infty}(a(y),b(y))|_{q=q_n}=(\hat{a}(y),\hat{b}(y))\ \text{in}\ C^1_{loc}([0,+\infty))\times C^1_{loc}([0,+\infty)),
\eeq
and
\beq
\left\{\arraycolsep=1.5pt
\begin{array}{lll}
d_S \hat{a}_{yy}+\hat{a}_y=0,\   \ &0<y<qL,\\[2mm]
d_I \hat{b}_{yy}+\hat{b}_y=0, \   \ &0<y<qL,\\[2mm]
d_S \hat{a}_y(0)+\hat{a}(0)=d_I \hat{b}_y(0)+\hat{b}(0)=0.\ &
\end{array}
\right.
\eeq
Direct calculation yields
\beq
\begin{cases}
\hat{a}(y)=K_Se^{-\f{y}{d_S}},\\
\hat{b}(y)=K_Ie^{-\f{y}{d_I}}.
\end{cases}
\eeq
Therefore
\beq
\label{lim an bn}
\lim\limits_{n\rightarrow+\infty}(a_n(y),b_n(y))=(K_Se^{-\f{y}{d_S}},K_Ie^{-\f{y}{d_I}})\ \text{in}\ C^1_{loc}([0,+\infty))\times C^1_{loc}([0,+\infty))
\eeq

\emph{Step 3.}
The values of $K_S$, $K_I$.

Due to the first equality of \eqref{10},
there holds
\beq
\Io\Lam(x)\dx=\I_0^{q_nL}\mu(L-\F{y}{q_n})a_n(y)\dy.
\eeq
Note from the first inequality of \eqref{a<= b<=} that
\beq
a_n(y)<\F{\|\Lam\|_1}{q^2}(1+\F{\be^*}{m\mu_*})+\F{\|\Lam\|_1e^{-\f{y}{d_S}}}{d_S(1-e^{-\f{qL}{d_S}})q}\big[\F{q}{\mu_*}+(1+\F{\be^*}{m\mu_*})L\big].
\eeq
Then, by a similar argument as the equality below (3.14) in \cite{cc},
we have that
\beq
\label{11}
\lim\limits_{n\rightarrow+\infty}\I_0^{q_nL}\mu(L-\F{y}{q_n})a_n(y)\dy
=\mu(L)\I_0^{+\infty}K_S e^{-\f{y}{d_S}}\dy.
\eeq
Hence $\|\Lam\|_1=\mu(L)\I_0^{+\infty}K_S e^{-\f{y}{d_S}}\dy$,
and moreover
$K_S=\F{\|\Lam\|_1}{d_S\mu(L)}$.

Let $v(x)=I(x)e^{-\f{q}{d_I}(1+\gam^*\f{d_I}{q^2})x}$.
Then from the $I$'s equation in \eqref{equilibrium}
we get that
\beq
\left\{\arraycolsep=1.5pt
\begin{array}{lll}
-d_I v_{xx}-q(1+2\F{\gam^*d_I}{q^2})v_x+(-\gam^*(1+\F{\gam^*d_I}{q^2})+\gam(x)-\F{\be(x)S(x)}{1+mI})v=0, \   \ &0<x<L,\\[2mm]
v_x=-\F{\gam^*}{q}v<0,\ &x=0,L.
\end{array}
\right.
\eeq
Take $v(x_0)=\min\limits_{x\in[0,L]}v(x)$.
Hence, $x_0\in(0,L]$,
since $v_x(0)<0$.
If $x_0\in(0,L)$,
then $v_x(x_0)=0,v_{xx}(x_0)\ge0$.
Therefore
$-\gam^*(1+\F{\gam^*d_I}{q^2})+\gam(x_0)-\F{\be(x_0)S(x_0)}{1+mI(x_0)}\ge0$,
which is impossible.
Hence $x_0=L$, $v(x)\ge v(L)$
and moreover
\beq
\F{I(x)}{I(L)}\ge e^{-\f{q}{d_I}(1+\gam^*\f{d_I}{q^2})(L-x)}.
\eeq

Take $I(x_1)=\max\limits_{x\in[0,L]}I(x)$.
Since $I_x(0)>0$,
$x_1\in(0,L]$.
If $x_1\in(0,L)$,
then $I_x(x_1)=0$.
By integrating both side of $I$'s equation in \eqref{equilibrium} from $0$ to $x_1$,
it indicates that
\beq
qI(x_1)=\I_0^{x_1}(\F{\be SI}{1+mI}-\gam I)\dx
\le I(x_1)\be^*\Io S\dx\le I(x_1)\be^*\F{\|\Lam\|_1}{\mu_*},
\eeq
which is impossible once $q>\F{\be^*\|\Lam\|_1}{\mu_*}$.
So, $x_1=L$ and $I(x)\le I(L)$ for all $x\in[0,L]$.
From the definition of $b_n$,
the inequality
\beq
e^{-\f{1}{d_I}(1+\gam^*\f{d_I}{q^2})y}\le\F{b_n(y)}{b_n(0)}\le1
\eeq
follows immediately.

Suppose that $K_I=0$,
i.e.,
$b_n(y)\rightarrow0$ in $C^1_{loc}([0,L))$ as $n\rightarrow+\infty$.
Define $\tilde{b}_n(y)=\F{b_n(y)}{b_n(0)}$.
Notice from the second equality of \eqref{10} that
\beq
\label{bb}
\I_0^{q_nL}\gam(L-\F{y}{q_n})\tilde{b}_n(y)\dy
=\I_0^{q_nL}\be(L-\F{y}{q_n})\F{a_n(y)}{\f{1}{q_n}+mb_n(y)}\tilde{b}_n(y)\dy.
\eeq
By the second equation of \eqref{a b's equ},
we get that
\beq
\left\{\arraycolsep=1.5pt
\begin{array}{ll}
d_I(\tilde{b}_n)_{yy}+(\tilde{b}_n)_y+\F{a_n\tilde{b}_n}{q_n(1+mq_nb_n)}\be(L-\F{y}{q_n})-\F{\tilde{b}_n}{q_n^2}\gam(L-\F{y}{q_n})=0, \   \ &0<y<q_nL,\\[2mm]
d_I (\tilde{b}_n)_y+\tilde{b}_n=0,\ &y=0,q_nL,\\[2mm]
\tilde{b}_n(y)\le1.&
\end{array}
\right.
\eeq
By the standard $L^p$-estimates of elliptic equations and Sobolev embedding theorem,
we obtain that
$\|\tilde{b}_n\|_{C^{1,\al}([0,M])}\le C$ for some $\al\in(0,1)$.
Passing to a subsequence again if necessary,
one can derive that
\beq
\label{lim tilde bn}
\lim\limits_{n\rightarrow+\infty}\tilde{b}_n:=b_{\infty}\ \text{in}\ C^1_{loc}([0,+\infty)),
\eeq
and
\beq
\left\{\arraycolsep=1.5pt
\begin{array}{lll}
d_I (b_{\infty})_{yy}+(b_{\infty})_y=0, \   \ &0<y<q_nL,\\[2mm]
d_I (b_{\infty})_y(0)+b_{\infty}(0)=0,\ &\\[2mm]
b_{\infty}(y)\le1,
\end{array}
\right.
\eeq
which implies that $b_{\infty}(y)=e^{-\f{y}{d_I}}$.
Similar to \eqref{11},
by some calculations,
one can have that
\beq
\label{13}
\lim\limits_{n\rightarrow+\infty}\I_0^{q_nL}\gam(L-\F{y}{q_n})\tilde{b}_n(y)\dy
=\gam(L)\I_0^{+\infty}e^{-\f{y}{d_S}}\dy
=d_I\gam(L).
\eeq
On the other hand,
once $q_n$ is large enough (e.g., $q_n>>\F{M}{L}$),
there holds $\be(L-\F{y}{q_n})>\F{\be(L)}{2}$ for $y\in[0,M]$.
Consequently,
it follows that
\beq
\I_0^{q_nL}\be(L-\F{y}{q_n})\F{a_n(y)}{\f{1}{q_n}+mb_n(y)}\tilde{b}_n(y)\dy\ge
\I_0^{M}\F12\be(L)\F{a_n(y)}{\f{1}{q_n}+mb_n(y)}\tilde{b}_n(y)\dy.
\eeq
This inequality, together with \eqref{lim an bn},
\eqref{lim tilde bn} and $q_n\rightarrow+\infty$,
yields that
\beq
\lim\limits_{n\rightarrow+\infty}\I_0^{q_nL}\be(L-\F{y}{q_n})\F{a_n(y)}{\f{1}{q_n}+mb_n(y)}\tilde{b}_n(y)\dy=+\infty,
\eeq
which contradicts to \eqref{bb} and \eqref{13}.
Therefore, it follows that $K_I>0$.

%Recall that
%\beq
%\I_0^{q_nL}\gam(L-\F{y}{q_n})b_n(y)\dy
%=\I_0^{q_nL}\be(L-\F{y}{q_n})\F{a_n(y)}{\f{1}{q_n}+mb_n(y)}b_n(y)\dy.
%\eeq
We next claim that
\beq
\label{12}
\bea{l}
\lim\limits_{n\rightarrow+\infty}
\Big(\I_0^{q_nL}\gam(L-\F{y}{q_n})b_n(y)\dy
-\I_0^{q_nL}\be(L-\F{y}{q_n})\F{a_n(y)}{\f{1}{q_n}+mb_n(y)}b_n(y)\dy\Big)\\
=\I_0^{+\infty}\Big(\be(L)\F{K_S}{m}e^{-\f{y}{d_S}}-\gam(L)K_Ie^{-\f{q}{d_I}}\Big)\dy.
\eea
\eeq
By \eqref{a<= b<=},
we have that there exists some positive integer $N>0$
such
that for all $n\ge N$,
\beq
\begin{cases}
a_n(y)\le\F{\|\Lam\|_1}{q_n^2}(1+\F{\be^*}{m\mu_*})
+\F{\|\Lam\|_1e^{-\f{y}{d_S}}}{d_S(1-e^{-\f{q_nL}{d_S}})}\big[\F{1}{\mu_*}+(1+\F{\be^*}{m\mu_*})\F{L}{q_n}\big],\\
b_n(y)\le\F{\|\Lam\|_1\be^*}{q_n^2m\mu_*}
+\F{\be^*\|\Lam\|_1e^{-\f{y}{d_I}}}{d_I m\mu_*(1-e^{-\f{q_nL}{d_I}})}\big[\F{1}{\gam_*}+\F{L}{q_n}\big].
\end{cases}
\eeq
For simplicity,
one can get the following estimates further for all $n\ge N$:
\beq
\begin{cases}
a_n(y)\le\F{c_1}{q_n^2}
+c_2e^{-\f{y}{d_S}},\\
b_n(y)\le\F{c_3}{q_n^2}
+c_4e^{-\f{y}{d_I}},
\end{cases}
\eeq
where $c_1,c_2,c_3,c_4$ are independent of $q_n$.
For any
$\eps>0$,
then there exists a $C(\eps)>0$ such that
\beq
\bea{l}
\I_{C(\eps)}^{q_nL}\Big(\be(L-\F{y}{q_n})\F{a_n(y)}{\f{1}{q_n}+mb_n(y)}b_n(y)+\gam(L-\F{y}{q_n})b_n(y)\Big)\dy
\le(\be^*+\gam^*)\I_{C(\eps)}^{q_nL}(q_na_n(y)b_n(y)+b_n(y))\dy\\
\le(\be^*+\gam^*)\Big[\I_0^{q_nL}(\F{c_1c_3}{q_n^3}+\F{c_3}{q_n^2})\dy+\I_{C(\eps)}^{+\infty}(\F{c_1c_4+c_2c_3}{q_n}e^{-\f{y}{d_S}}+c_2c_4e^{-\f{2y}{d_I}})\dy\Big]<\F{\eps}{3},
\eea
\eeq
and
\beq
\I_{C(\eps)}^{q_nL}\Big(\be(L)\F{K_S}{m}e^{-\f{y}{d_S}}-\gam(L)K_Ie^{-\f{q}{d_I}}\Big)\dy
\le(\be^*+\gam^*)\I_{C(\eps)}^{q_nL}\Big(\F{K_S}{m}e^{-\f{y}{d_S}}+K_Ie^{-\f{q}{d_I}}\Big)\dy<\F{\eps}{3}.
\eeq
Since
$\lim\limits_{n\rightarrow+\infty}\F{a_n(y)b_n(y)}{\f{1}{q_n}+mb_n(y)}=\F{K_S}{m}e^{-\f{y}{d_S}}$
and
$\lim\limits_{n\rightarrow+\infty}b_n(y)=K_Ie^{-\f{y}{d_I}}$ in $C^1_{loc}([0,+\infty))$,
it then follows that
\beq
\begin{cases}
\I_0^{C(\eps)}\Big|\be(L-\F{y}{q_n})\F{a_n(y)b_n(y)}{\f{1}{q_n}+mb_n(y)}-\be(L)\F{K_S}{m}e^{-\f{y}{d_S}}\Big|\dy<\F{\eps}{6},\\
\I_0^{C(\eps)}\Big|\gam(L-\F{y}{q_n})b_n(y)-\gam(L)K_Ie^{-\f{q}{d_I}}\Big|\dy<\F{\eps}{6},
\end{cases}
\eeq
for large $n$.
Therefore,
it is obvious that there exists a $N>0$ such that for all $n>N$,
\beq
\bea{l}
\Big|\Big(\I_0^{q_nL}\gam(L-\F{y}{q_n})b_n(y)\dy
-\I_0^{q_nL}\be(L-\F{y}{q_n})\F{a_n(y)}{\f{1}{q_n}+mb_n(y)}b_n(y)\dy\Big)\\
-\I_0^{+\infty}\Big(\be(L)\F{K_S}{m}e^{-\f{y}{d_S}}-\gam(L)K_Ie^{-\f{q}{d_I}}\Big)\dy\Big|\\
\le\I_0^{C(\eps)}\Big|\be(L-\F{y}{q_n})\F{a_n(y)b_n(y)}{\f{1}{q_n}+mb_n(y)}-\be(L)\F{K_S}{m}e^{-\f{y}{d_S}}\Big|\dy
+\I_0^{C(\eps)}\Big|\gam(L-\F{y}{q_n})b_n(y)-\gam(L)K_Ie^{-\f{q}{d_I}}\Big|\dy\\
+\I^{q_nL}_{C(\eps)}\Big|\be(L-\F{y}{q_n})\F{a_n(y)}{\f{1}{q_n}+mb_n(y)}b_n(y)+\gam(L-\F{y}{q_n})b_n(y)\Big|\dy+\I^{q_nL}_{C(\eps)}\Big|\be(L)\F{K_S}{m}e^{-\f{y}{d_S}}-\gam(L)K_Ie^{-\f{q}{d_I}}\Big|\dy<\eps,
\eea
\eeq
which proves our claim \eqref{12}.
Hence,
\beq
K_I=\F{\be(L)K_S\i_0^{+\infty}e^{-\f{y}{d_S}}\dy}{m\gam(L)\i_0^{+\infty}e^{-\f{y}{d_I}}\dy}
=\F{\be(L)\|\Lam\|_1}{m\gam(L)d_I\mu(L)}.
\eeq
We therefore conclude this theorem.
\proofend

\subsection{Concentration behaviors of EE as $d_S\rightarrow0$}

This subsection is devoted to proving Theorem \ref{ee ds-0},
i.e.,
we shall give the asymptotic profiles of EE for small $d_S$.

\four
For the sake of clarity, the proof will be divided into several steps.

\emph{Step1.}
Introduce
$W(y)=d_S S(L-d_Sy)$.
It can be deduced from \eqref{s<=} that for sufficiently small $d_S$
\beq
\begin{cases}
W(y)\le\F{d_S\|\Lam\|_1}{q}(1+\F{\be^*}{m\mu_*})
+\F{\|\Lam\|_1e^{-qy}}{(1-e^{-\f{qL}{d_S}})}\big[\F{q}{\mu_*}+(1+\F{\be^*}{m\mu_*})L\big],~\forall~0\le y\le\F{L}{d_S},\\
W_{yy}+qW_y+d^2_S\Lam(L-d_Sy)-d_S\mu(L-d_Sy)W-\F{d_S\be WI}{1+mI}+d^2_S\gam=0,~\forall~0\le y\le\F{L}{d_S},\\
W_y+qW=0,~y=0,\F{L}{d_S}.
\end{cases}
\eeq
Note that
\beq
I(L-d_Sy)\le\F{\|\Lam\|_1\be^*}{qm\mu_*}
+\F{\be^*\|\Lam\|_1e^{-\f{d_Sqy}{d_I}}}{d_I m\mu_*(1-e^{-\f{qL}{d_I}})}\big[\F{q}{\gam_*}+L\big].
\eeq
Consequently, for any fixed constant $M>0$,
$W$ and $I$ are uniformly bounded in
$[0,M]$.
Applying the standard elliptic regularity and Sobolev embedding
theorem,
we get $\{d_{S_k}\}$ satisfying $\lim\limits_{k\rightarrow+\infty}d_{S_k}=0$
by passing to a subsequence if necessary,
and have that
\beq
\lim\limits_{k\rightarrow+\infty}W_k(y):=\lim\limits_{k\rightarrow+\infty}W(y)|_{d_S=d_{S_k}}=\hat{W}(y)\ \text{in}\ C^1_{loc}([0,+\infty)).
\eeq
Here $\hat{W}$ satisfies
\beq
\begin{cases}
\hat{W}_{yy}+q\hat{W}_y=0,~0\le y<+\infty,\\
\hat{W}_y(0)+q\hat{W}(0)=0.
\end{cases}
\eeq
Then $\hat{W}(y)=c_Se^{-qy}$ for some constant $c_S\ge0$.
Note that $\Io\Lam(x)\dx=\I_0^{\f{L}{d_S}}\mu(L-d_Sy)W(y)\dy$.
Similar to \eqref{11},
it can be proved that $\Io\Lam(x)\dx=\I_0^{+\infty}\mu(L)c_Se^{-qy}\dy$,
which implies $c_S=\F{q\|\Lam\|_1}{\mu(L)}$.
Hence,  we complete the proof of statement (1).

\emph{Step 2.}
We denote
$\tilde{I}(x)=I(x)e^{-\f{qx}{d_I}}$. Then it follows
from $I's$ equation in \eqref{equilibrium} that
\beq
\left\{\arraycolsep=1.5pt
\begin{array}{ll}
d_I\tilde{I}_{xx}+q\tilde{I}_x+\F{\be S\tilde{I}}{1+mI}-\gam(x)\tilde{I}=0, \   \ &0<x<L,\\[2mm]
\tilde{I}_x=0,\ &x=0,L.
\end{array}
\right.
\eeq
Owing to
\beq
\Io\big|\F{\be S}{1+mI}-\gam(x)\big|\dx\le\be^*\Io S(x)\dx+\gam^*L\le\F{\be^*\|\Lam\|_1}{\mu_*}+\gam^*L,
\eeq
the Harnack inequality \cite[Lemma 2.2]{l} yields that
\beq
\max\limits_{x\in[0,L]}\tilde{I}(x)\le C\min\limits_{x\in[0,L]}\tilde{I}(x)\ \text{for any}\ d_S>0.
\eeq
Hence there also holds
\beq
\label{harnack}
\max\limits_{x\in[0,L]}I(x)\le C\min\limits_{x\in[0,L]}I(x)\ \text{for any}\ d_S>0.
\eeq
Note from \eqref{i<=} that $I(x)$ is uniformly bounded on $[0,L]$ for all $d_S>0$. Then
\beq
\Io\big|\F{\be SI}{1+mI}-\gam(x)I\big(\dx\le\Io\big|\F{\be S}{m}+\gam(x)I\big)\dx<+\infty.
\eeq
By $L^1$ regularity and Sobolev embedding theorem,
we may assume that,
passing to a
subsequence if necessary, there exists a nonnegative function $I^{\infty}(x)\in C^{0,\al}([0,L])$ for some $\al\in(0,1)$
such that $I(x)\rightarrow I^{\infty}(x)$ in $C^{0,\al}([0,L])$.
Moreover,
\eqref{harnack} yields that
$I^{\infty}>0$ or $I^{\infty}\equiv0$ on $[0,L]$.

\emph{Step 3.}
By \eqref{s<=},
we can
see that $S(x)$ is locally uniformly bounded in $[0,L)$ for all $d_S>0$.
Then for any fixed small $\eps>0$,
passing to a subsequence of $d_S\rightarrow0$ if necessary,
there
exists $S^{\infty}(x)\in L^2((0,L-\eps))$ such that $S(x)\rightarrow S^{\infty}(x)\ge0$ weakly in $L^2((0,L-\eps))$.
Note that $\Io S(x)\dx\le\F{\|\Lam\|_1}{\mu_*}$. By a diagonal argument,
up to a further subsequence,we arrive at
\beq
\Io S^{\infty}(x)\dx=\lim\limits_{\eps\rightarrow0}\I_0^{L-\eps}S^{\infty}(x)\dx
=\lim\limits_{\eps\rightarrow0}\lim\limits_{d_S\rightarrow0}\I_0^{L-\eps}S(x)\dx
\le\F{\|\Lam\|_1}{\mu_*}.
\eeq
Then
$S^{\infty}\in L^1((0,L))$.

Choosing any test function $\phi\in C_0^2([0,L-\eps])$,
multiplying $S's$ equation by $\phi$ and integrating by parts over $(0,L-\eps)$,
it follows that
\beq
d_S\I_0^{L-\eps}S\phi_{xx}\dx+q\I_0^{L-\eps}S\phi_x\dx+\I_0^{L-\eps}\Lam\phi\dx-\I_0^{L-\eps}\mu S\phi\dx-\I_0^{L-\eps}\F{\be SI}{1+mI}\phi\dx+\I_0^{L-\eps}\gam I\phi\dx=0.
\eeq
Hence, by letting $d_S\rightarrow0$,
because of the arbitrariness of $\eps$,
one can have that
$S^{\infty}$ is a weak
solution and
then a classical solution of the first equation in
\eqref{sinfty iinfty}.

\emph{Step 4.}
Let $S(x_{d_S})=\min\limits_{x\in[0,L]}S(x)$.
Then by $S_x(0)>0$,
$x_{d_S}\in[0,L)$.
We claim that $\lim\limits_{d_S\rightarrow0}S(x_{d_S})=0$.
If it is false,
then there exist $\eps_0>0$ and $0<d_{S_n}<\F{1}{n}$ such that $S(x_{d_{S_n}})>\eps_0$ for all $n\in\mathrm{N}$.
Set $w(x)=e^{-\f{q}{d_{S_n}}(1-c\f{d_{S_n}}{q^2})}S(x)$.
It follows that
\beq
\left\{\arraycolsep=1.5pt
\begin{array}{ll}
d_{S_n} w_{xx}+q(1-2c\F{d_{S_n}}{q^2})w_x+\Big(-c(1-c\F{d_{S_n}}{q^2})-\mu(x)-\F{\be(x)S(x)I(x)}{1+mI}+\F{\Lam(x)+\gam(x)I}{S(x)}\Big)w=0, \   \ &0<x<L,\\[2mm]
w_x=\F{c}{q}w>0,\ &x=0,L.
\end{array}
\right.
\eeq
Let $w(\hat{x})=\max\limits_{x\in[0,L]}w(x)$.
Clearly, $\hat{x}\in(0,L]$.
If $\hat{x}\in(0,L)$,
it can be easily checked that $w_x(\hat{x})=0$, $w_{xx}(\hat{x})\le0$.
Therefore,
\beq
-c(1-c\F{d_{S_n}}{q^2})-\mu(\hat{x})-\F{\be(\hat{x})S(\hat{x})I(\hat{x})}{1+mI(\hat{x})}+\F{\Lam(\hat{x})+\gam(\hat{x})I(\hat{x})}{S(\hat{x})}\ge0.
\eeq
However,
\beq
\bea{l}
-c(1-c\F{d_{S_n}}{q^2})-\mu(\hat{x})-\F{\be(\hat{x})S(\hat{x})I(\hat{x})}{1+mI(\hat{x})}+\F{\Lam(\hat{x})+\gam(\hat{x})I(\hat{x})}{S(\hat{x})}
\le-c(1-c\F{d_{S_n}}{q^2})+\F{\Lam(\hat{x})+\gam(\hat{x})I(\hat{x})}{\eps_0}\\
\le-c(1-c\F{d_{S_n}}{q^2})+\F{\Lam^*+\gam^*I^*}{\eps_0}<0\ \text{as}\ n\rightarrow+\infty\ \text{and}\ c>\F{\Lam^*+\gam^*I^*}{\eps_0}.
\eea
\eeq
A contradiction then occurs.
Hence $\hat{x}=L$.
Then $w(x)\le w(L)$,
which indicates that
\beq
S(x)\le S(L)e^{-\f{q}{d_{S_n}}(1-c\f{d_{S_n}}{q^2})(L-x)}\ \text{for all}\ x\in[0,L].
\eeq
By the estimate of $S(L)$ in \eqref{a},
it holds
\beq
S(x)\le \F{1}{d_{S_n}(1-e^{-\f{qL}{d_{S_n}}})}\big[\F{\|q\Lam\|_1}{\mu_*}+\|\Lam\|_1(1+\F{\be^*}{m\mu_*})L\big]e^{-\f{q}{d_{S_n}}(1-c\f{d_{S_n}}{q^2})(L-x)}.
\eeq
In conclusion,we have
$S(x)\rightarrow0$ locally uniformly in $[0,L)$ as $n\rightarrow+\infty$,
which contradicts to $S(x)\ge\eps_0$.
Thus $\lim\limits_{d_S\rightarrow0}S(x_{d_S})=0$.

\emph{Step 5.}
$x_{d_S}=0$.

Suppose that $x_{d_S}\in(0,L)$. Then $S_x(x_{d_S})=0$, $S_{xx}(x_{d_S})\ge0$.
Now by using $S's$ equation,
we
obtain
\beq
\Lam(x_{d_S})-\mu(x_{d_S})S(x_{d_S})+\gam(x_{d_S})I(x_{d_S})-\F{\be(x_{d_S})S(x_{d_S})I(x_{d_S})}{1+mI(x_{d_S})}\le0.
\eeq
This inequality implies obviously that
\beq
\Lam_*\le\Lam(x_{d_S})\le\mu(x_{d_S})S(x_{d_S})+\F{\be(x_{d_S})S(x_{d_S})I(x_{d_S})}{1+mI(x_{d_S})}
\le(\mu^*+\F{\be^*}{m})S(x_{d_S})\rightarrow0\ \text{as}\ d_S\rightarrow0,
\eeq
which leads to a contradiction.
So, $x_{d_S}=0$.
Therefore $S^{\infty}(0)=\lim\limits_{d_S\rightarrow0}S(0)=\lim\limits_{d_S\rightarrow0}S(x_{d_S})=0$.
In conclusion, $S^{\infty}$ is a solution of
\beq
\label{b}
\left\{\arraycolsep=1.5pt
\begin{array}{ll}
-q S^{\infty}_x+\Lam(x)-\mu(x)S^{\infty}-\F{\be(x)S^{\infty}I^{\infty}}{1+mI^{\infty}}+\gam(x)I^{\infty}=0,\   \ &0<x<L,\\[2mm]
S^{\infty}(0)=0,\ &x=0,L.
\end{array}
\right.
\eeq

\emph{Step 6.}
The equation of $I^{\infty}$.

By the similar arguments as above,
we have that
$I^{\infty}(x)$ is a classical solution of
\beq
d_I I^{\infty}_{xx}-q I_x+\F{\be S^{\infty}I^{\infty}}{1+mI^{\infty}}-\gam(x)I^{\infty}=0.
\eeq
Since for any given $\eps\in(0,L)$, $I(x)\in L^{\infty}([0,L])$,
$S(x)\in L^{\infty}([0,L-\eps])$ for all $d_S>0$,
by means of $I's$ equation and the standard elliptic regularity and Sobolev embedding theorem,
passing to a
subsequence if necessary
there holds
$I(x)\rightarrow I^{\infty}(x)$ in $C^1([0,L-\eps])$.
In particular,
$I(0)\rightarrow I^{\infty}(0)$,
$I_x(0)\rightarrow I_x^{\infty}(0)$.
Thus,
$d_I I^{\infty}_x(0)-qI^{\infty}(0)=0$.
It then follows that $I^{\infty}$ is a classical solution of \beq
\label{iinfty}
\left\{\arraycolsep=1.5pt
\begin{array}{ll}
d_I I^{\infty}_{xx}-q I_x+\F{\be S^{\infty}I^{\infty}}{1+mI^{\infty}}-\gam(x)I^{\infty}=0, \   \ &0<x<L,\\[2mm]
d_I I^{\infty}_x(0)-qI^{\infty}(0)=0.&
\end{array}
\right.
\eeq

Now we can use the similar arguments as to the Line $5-13$ on page 20 in \cite{cc} to deduce that
\beq
S(x)\rightarrow S^{\infty}+N_S\del_{L}(x)\ \text{in}\ L^1((0,L)),
\eeq
and in fact $N_S=\lim\limits_{\eps\rightarrow0}\lim\limits_{d_S\rightarrow0}\I_{L-\eps}^{L}S(x)\dx$,
where $N_S$ is given in \eqref{ns}.

\emph{Step 7.}
The boundary condition at downstream end.

Integrating the $I's$ equation in $(L-\eps,L)$ yields
\beq
\bea{l}
d_I I_x(L-\eps)-qI(L-\eps)=\I_{L-\eps}^{L}\F{\be SI}{1+mI}\dx-\I_{L-\eps}^{L}\gam I\dx\\
=\F{\be(x_{\eps})I(x_{\eps})}{1+mI(x_{\eps})}\I_{L-\eps}^{L}S(x)\dx-\I_{L-\eps}^L\gam(x)I(x)\dx\ \text{for some}\ x_{\eps}\in(L-\eps,L).
\eea
\eeq
By letting $d_S\rightarrow0$ first and then taking $\eps\rightarrow0$,
it is not hard to deduce that
\beq
d_I I_x^{\infty}(L)-qI^{\infty}(L)=N_S\F{\be(L)I^{\infty}(L)}{1+mI^{\infty}(L)}.
\eeq
Adding the $S's$ equation and the $I's$ equation up,
integrating them in $(0,L)$,
it
follows that \beq
d_I I^{\infty}_x(L)-qI^{\infty}-qS^{\infty}(L)+\Io(\Lam(x)-\mu(x)S^{\infty}(x))\dx
=N_S\F{\be(L)I^{\infty}(L)}{1+mI^{\infty}(L)}-qS^{\infty}(L)+N_S\mu(L),
\eeq
which implies that
\beq
\label{sinfty(L)}
qS^{\infty}(L)=N_S(\mu(L)+\F{\be(L)I^{\infty}(L)}{1+mI^{\infty}(L)}).
\eeq

\emph{Step 8: $S^{\infty}(x)>0$ in $(0,L]$, $I^{\infty}(x)>0$ on $[0,L]$.}

Suppose that there exists a $\hat{x}\in(0,L)$ such that $S^{\infty}(\hat{x})=0$,
then $S^{\infty}_x(\hat{x})\le0$.
From \eqref{b} we can see that
$\Lam(\hat{x})+\gam(\hat{x})I^{\infty}(\hat{x})=0$,
which is impossible.
If $S^{\infty}(L)=0$,
then this contradicts to statement (1).

Recall that $I^{\infty}>0$ or $I^{\infty}\equiv0$.
Suppose that $I^{\infty}\equiv0$ on $[0,L]$.
Then there exist $d_{S_K}$, $S_k(x)$, $I_k(x)$ such that $d_{S_k}\rightarrow0$ and $I_k\rightarrow0$ as $k\rightarrow+\infty$.
Set $\tilde{I}_k=\F{I_k}{\|I_k\|_{\infty}}$.
We can get that
\beq
\left\{\arraycolsep=1.5pt
\begin{array}{ll}
d_I (\tilde{I}_k)_{xx}-q (\tilde{I}_k)_x+(\be(x)\F{S_k}{1+mI_k}-\gam(x))\tilde{I}_k=0, \   \ &0<x<L,\\[2mm]
d_I (\tilde{I}_k)_x-q \tilde{I}_k=0,\ &x=0,L,\\[2mm]
\|\tilde{I}_k\|_{\infty}=1.&
\end{array}
\right.
\eeq
Since $\be(x)\F{S_k}{1+mI_k}-\gam(x)\in L^1((0,L))$,
then after passing to a further subsequence
if necessary,
one can assume that $\tilde{I}_k\rightarrow\tilde{I}^{\infty}\ge0$ in $C^{\al}([0,L])$ and $\|\tilde{I}^{\infty}\|_{\infty}=1$ as $k\rightarrow+\infty$.
Hence, $\tilde{I}^{\infty}$ is a solution of
\beq
\left\{\arraycolsep=1.5pt
\begin{array}{ll}
d_I (\tilde{I}^{\infty})_{xx}-q (\tilde{I}^{\infty})_x+(\be\tilde{S}^{\infty}-\gam(x))\tilde{I}^{\infty}=0, \   \ &0<x<L,\\[2mm]
d_I (\tilde{I}^{\infty})_x(0)-q \tilde{I}^{\infty}(0)=0,\ &\\[2mm]
d_I (\tilde{I}^{\infty})_x(L)=\tilde{I}^{\infty}(L)(q+N_S\be(L)).&
\end{array}
\right.
\eeq
Actually $\tilde{S}^{\infty}$ is the solution of \eqref{s infty}.
This indicates that the eigenvalue problem \eqref{aux eigen} has the principal eigenvalue $0$,
a contradiction with the assumption \eqref{assum}.
Hence $S^{\infty}>0,I^{\infty}>0$.
\proofend

\subsection{Concentration behaviors of EE as $d_I\rightarrow0$}

\five
We split our proof into three steps.

\emph{Step 1. Convergence of $I$.}

According to \eqref{s<=},
we get that for all $x\in[0,L]$
\beq
\label{s<=A}
S(x)\le\F{\|\Lam\|_1}{q}(1+\F{\be^*}{m\mu_*})
+\F{\|\Lam\|_1}{d_S(1-e^{-\f{qL}{d_S}})}\big[\F{q}{\mu_*}+(1+\F{\be^*}{m\mu_*})L\big]:=A.
\eeq
Let
$Y(x)=I(x)e^{-\f{q}{d_I}(1-c\f{d_I}{q^2})x}$.
By $I$'s equation,
there holds
\beq
\left\{\arraycolsep=1.5pt
\begin{array}{ll}
-d_I Y_{xx}-q(1-2c\F{d_I}{q^2})Y_x+\Big(c(1-c\F{d_I}{q^2})+\gam(x)-\F{\be(x)S(x)}{1+mI}\Big)Y=0, \   \ &0<x<L,\\[2mm]
Y_x=\F{c}{q}Y>0,\ &x=0,L.
\end{array}
\right.
\eeq
Take $x_1\in[0,L]$ such that $Y(x_1)=\max\limits_{x\in[0,L]}Y(x)$.
Since $Y_x(0)>0$,
then $x_1>0$.
By noting the equation of $Y$ above,
if $x_1\in(0,L)$,
it follows that
$Y_x(x_1)=0,Y_{xx}(x_1)\le0$ and
\beq
c(1-c\F{d_I}{q^2})+\gam(x_1)-\F{\be(x_1)S(x_1)}{1+mI(x_1)}\le0.
\eeq
On the other hand,
together with \eqref{s<=A} and sufficiently small $d_I$,
we reach that
\beq
c(1-c\F{d_I}{q^2})+\gam(x_1)-\F{\be(x_1)S(x_1)}{1+mI(x_1)}
\ge c(1-c\F{d_I}{q^2})+\gam_*-\be^*A
\ge\F{c}{2}+\gam_*-\be^*A>0,
\eeq
provided that $c>\max\{0,2(\be^*A-\gam_*)\}$.
A contradiction then occurs.
Hence $x_1=L$.
Therefore, $Y(x)\le Y(L)$,
i.e.,
$I(x)\le I(L)e^{-\f{q}{d_I}(1-c\f{d_I}{q^2})(L-x)}$.
Combined with the estimate of $I(L)$ in \eqref{ix<= iL<=},
one can obtain that
\beq
\label{14}
\bea{l}
I(x)\le\F{\be^*\|\Lam\|_1}{d_Im\mu_*(1-e^{-\f{qL}{d_I}})}\big[\F{q}{\gam_*}+(L-\F{d_I}{q}(1-e^{-\f{qL}{d_I}}))\big]e^{-\f{q}{d_I}(1-c\f{d_I}{q^2})(L-x)}\\
\le\F{\be^*\|\Lam\|_1}{d_Im\mu_*(1-e^{-\f{qL}{d_I}})}\big[\F{q}{\gam_*}+L\big]e^{-\f{q}{d_I}(1-c\f{d_I}{q^2})(L-x)},\ \text{for all}\ x\in[0,L].
\eea
\eeq
Hence, $I(x)\rightarrow0$ locally uniformly
in $[0,L)$ as $d_I\rightarrow0$.

\emph{Step 2 Convergence of $\Io I(x)\dx$.}

Since $\|\F{\be(x)S(x)}{1+mI(x)}\|_{\infty}\le\|\be(x)S(x)\|_{\infty}\le\be^*A$,
and $\F{\be(x)S(x)I(x)}{1+mI(x)}<\be^*AI(x)\rightarrow0$ for $x\in[0,L)$ as $d_I\rightarrow0$,
integrating $I$'s equation from $0$ to $L$ yields that
\beq
\label{convergence of int I}
\gam_*\lim\limits_{d_I\rightarrow0}
\Io I(x)\dx\le\lim\limits_{d_I\rightarrow0}\Io\F{\be(x)S(x)I(x)}{1+mI(x)}\dx=0,
\eeq
where we have applied the Lebesgue-dominated convergence theorem.

\emph{Step 3 Convergence of $S$.}

By \eqref{14},
we have that
\beq
\Io I(x)\dx\le\F{\be^*\|\Lam\|_1}{m\mu_*(1-e^{-\f{qL}{d_I}})q(1-c\f{d_I}{q^2})}\big[\F{q}{\gam_*}+L\big](1-e^{-\f{q}{d_I}(1-c\f{d_I}{q^2})L}).
\eeq

Combining  \eqref{s<=A}, the inequality above with the
elliptic $L^1$-theory in \cite{bs},  \eqref{equilibrium} yields
$S(x)\in W^{1,p}((0,L))$ for any $p\ge1$.

By means of the Sobolev embedding theorem,
after passing to a subsequence of $d_I\rightarrow0$,
we may assume that
\beq
S(x)\rightarrow\hat{S}(x)\ge0\ \text{in}\ C([0,L]),\ \text{when}\ d_I\rightarrow0.
\eeq
Taking any test function $\phi\in C_0^2((0,L))$,
multiplying $S's$ equation by $\phi$ and integrating by parts over $(0,L)$,
it follows that
\beq
d_S\I_0^{L}S\phi_{xx}\dx+q\I_0^{L}S\phi_x\dx+\I_0^{L}\Lam\phi\dx-\I_0^{L}\mu S\phi\dx-\I_0^{L}\F{\be SI}{1+mI}\phi\dx+\I_0^{L}\gam I\phi\dx=0.
\eeq
By \eqref{convergence of int I},
one can obtain that
\beq
\lim\limits_{d_I\rightarrow0}\Io|\gam(x)-\F{\be(x)S(x)}{1+mI(x)}|I\phi\dx\le\lim\limits_{d_I\rightarrow0}(\gam^*+\be^*A)\|\phi\|_{\infty}\Io I(x)\dx=0.
\eeq
It is evident that $\hat{S}$ satisfies
\beq
d_S\hat{S}_{xx}-q\hat{S}_x+\Lam(x)-\mu(x)\hat{S}=0,\ \forall\ x\in(0,L).
\eeq

Note that for any given $\eps\in(0,L)$, $I(x)\in L^{\infty}([0,L-\eps])$,
$S(x)\in L^{\infty}([0,L])$ for all $d_S>0$.
By means of $S$'s equation and the standard elliptic regularity and Sobolev embedding theorem,
passing to a
subsequence if necessary,
there holds
$S(x)\rightarrow \hat{S}(x)$ in $C^1([0,L-\eps])$
when $d_I\rightarrow0$.
In particular,
$0=d_S S_x(0)-qS(0)\rightarrow d_S\hat{S}_x(0)-q\hat{S}(0)$.
It then follows that $\hat{S}$ is a classical solution of
\beq
\left\{\arraycolsep=1.5pt
\begin{array}{ll}
d_S\hat{S}_{xx}-q\hat{S}_x+\Lam(x)-\mu(x)\hat{S}=0, \   \ &0<x<L,\\[2mm]
d_S \hat{S}_x(0)-q\hat{S}(0)=0.&
\end{array}
\right.
\eeq

Finally,
we check the boundary condition at downstream end.
Integrating the $S's$ equation in $(L-\eps,L)$ yields
\beq
\bea{l}
d_S S_x(L-\eps)-qS(L-\eps)=\I_{L-\eps}^{L}\Lam(x)\dx-\I_{L-\eps}^{L}\mu(x)S(x)\dx-\I_{L-\eps}^{L}\F{\be SI}{1+mI}\dx+\I_{L-\eps}^{L}\gam I\dx\\
=\I_{L-\eps}^{L}\Lam(x)\dx-\mu(x_{\eps})\I_{L-\eps}^{L}S(x)\dx-\F{\be(x_{\eps})I(x_{\eps})}{1+mI(x_{\eps})}\I_{L-\eps}^{L}S(x)\dx+\I_{L-\eps}^L\gam(x)I(x)\dx,
\eea
\eeq
for some $x_{\eps}\in(L-\eps,L)$.
By letting $d_I\rightarrow0$ first and then taking $\eps\rightarrow0$,
it is not hard to deduce that
\beq
d_S \hat{S}_x(L)-q\hat{S}(L)=0.
\eeq
Hence $\hat{S}$ is a classical solution of \eqref{hattilde S},
Furthermore we must have $\hat{S}(x)=\s(x)$ and $\lim\limits_{d_I\rightarrow0}S(x)=\s(x)$ uniformly on $[0,L]$.
The proof of Theorem \ref{ee di-0} is then finished.
\proofend

\subsection{Concentration behaviors of EE as $d_S\rightarrow+\infty$}

\six

By \eqref{s<=} and \eqref{i<=},
we get that $S(x)$ and $I(x)$ are bounded uniformly for $x\in[0,L]$ as $d_S\rightarrow+\infty$.
Then by passing to a subsequence of $d_S\rightarrow+\infty$ if necessary,
the standard elliptic regularity and Sobolev imbedding
theorem guarantee that
$(S(x),I(x))\rightarrow(S^{\infty}(x),I^{\infty}(x))$ in $C^1([0,L])\times C^1([0,L])$.
Integrating the $S$'s equation in \eqref{equilibrium}
from $0$ to $x$ and dividing both side by $d_S$ yield that
\beq
S_x=\F{q}{d_S}S(x)+\I_0^x\big(-\Lam(x)+\mu(x)S(x)+\F{\be(x)S(x)I(x)}{1+mI(x)}-\gam(x)I(x)\big)\dx\rightarrow0\ \text{as}\ d_S\rightarrow+\infty,
\eeq
which implies that $S^{\infty}(x)$ is a constant.
Then by noting the first identity in \eqref{10},
it follows
that $S^{\infty}=\F{\|\Lam\|_1}{\|\mu\|_1}$.
By applying the Harnack inequality in $I$'s equation,
there holds $\max\limits_{x\in[0,L]}I(x)\le C\min\limits_{x\in[0,L]}I(x)$.
Hence $I^{\infty}\equiv0$ or $I^{\infty}>0$ on $[0,L]$.

Suppose that $I^{\infty}\equiv0$.
Let $\tilde{I}:=\F{I}{\|I\|_{\infty}}$,
and then $\|\tilde{I}\|_{\infty}=1$.
Therefore $\tilde{I}$ satisfies
\beq
\left\{\arraycolsep=1.5pt
\begin{array}{ll}
d_I \tilde{I}_{xx}-q \tilde{I}_x+(\beta(x)\F{S}{1+mI}-\gam(x))\tilde{I}=0, \   \ &0<x<L,\\[2mm]
d_I \tilde{I}_x-q \tilde{I}=0,\ &x=0,L.
\end{array}
\right.
\eeq
By passing to a subsequence of $d_S\rightarrow+\infty$ if necessary,
we may assume from  the elliptic regularity and Sobolev embedding
theorem  that
$\tilde{I}\rightarrow\hat{I}$ in $C^1([0,L])$ as $d_S\rightarrow+\infty$,
where $\hat{I}\ge0$,
$\|\hat{I}\|_{\infty}=1$ and satisfies
\beq
\left\{\arraycolsep=1.5pt
\begin{array}{ll}
d_I \hat{I}_{xx}-q \hat{I}_x+(\beta(x)S^{\infty}(x)-\gam(x))\hat{I}=0, \   \ &0<x<L,\\[2mm]
d_I \hat{I}_x-q \hat{I}=0,\ &x=0,L.
\end{array}
\right.
\eeq
It then follows from the Harnack inequality that $\hat{I}>0$ on $[0,L]$.
This indicates that $\lam_1(d_I,q,\gam-\be S^{\infty})=0$,
and furthermore $\rr_0^*=1$ owing to the
statement (i) of Theorem \ref{properties of R0}.
A contradiction then occurs to our assumption $\rr_0^*>1$
(i.e.,$\lam_1(d_I,q,\gam-\be S^{\infty})$ should be negative). Thus, $I^{\infty}>0$
on $[0,L]$ and satisfies
\beq
\label{15}
\left\{\arraycolsep=1.5pt
\begin{array}{ll}
d_I I^{\infty}_{xx}-q I^{\infty}_x+\beta(x)\F{S^{\infty}I^{\infty}}{1+mI^{\infty}}-\gam(x)I^{\infty}=0, \   \ &0<x<L,\\[2mm]
d_I I^{\infty}_x-q I^{\infty}=0,\ &x=0,L.
\end{array}
\right.
\eeq

Finally,
we show the uniqueness of $I^{\infty}$ as follows.
Let $\mathcal{I}:=I^{\infty}e^{-\f{qx}{d_I}}$.
Then \eqref{15} can be rewritten as
\beq
\label{16}
\left\{\arraycolsep=1.5pt
\begin{array}{ll}
d_I \mathcal{I}^{\infty}_{xx}+q \mathcal{I}^{\infty}_x+\beta(x)\F{S^{\infty}\mathcal{I}^{\infty}}{1+me^{\f{qx}{d_I}}\mathcal{I}^{\infty}}-\gam(x)\mathcal{I}^{\infty}=0, \   \ &0<x<L,\\[2mm]
\mathcal{I}^{\infty}_x=0,\ &x=0,L.
\end{array}
\right.
\eeq
Suppose on the contrary that
there are two different positive solution of \eqref{16},
denoted by $\mathcal{I}_i$,
$i=1,2$.
Since the equation above has arbitrary large super-solutions (e.g., any constant $C$ large enough) and arbitrary small positive sub-solutions,
e.g., $\eps\psi$, where $\eps>0$ small and $\psi>0$ is an eigenfunction of
\beq
\left\{\arraycolsep=1.5pt
\begin{array}{ll}
-d_I\psi_{xx}-q\psi_x+(\gam(x)-S^{\infty}(x)\be(x))\psi=\lam\psi,\   \ &0<x<L,\\[2mm]
\psi_x=0,\ &x=0,L,
\end{array}
\right.
\eeq
then we can choose
$\eps$ and $C$ such that $\eps\psi\le \mathcal{I}_1,\mathcal{I}_2\le C$.
Hence, the maximal solution and the minimal
solution,
denoted by $\mathcal{I}_m$ and $\mathcal{I}_M$,
also satisfy $\mathcal{I}_m\le\mathcal{I}_1,\mathcal{I}_2\le \mathcal{I}_M$ in $(0,L)$.
Since $\mathcal{I}_1\not\equiv\mathcal{I}_2$,
we have $\mathcal{I}_M\ge,\not\equiv\mathcal{I}_m$. Multiplying the equation of $\mathcal{I}_m$ by $\mathcal{I}_Me^{\f{qx}{d_I}}$,
the equation of $\mathcal{I}_M$ by $\mathcal{I}_me^{\f{qx}{d_I}}$,
subtracting and integrating the result in $(0,L)$,
we see that
\beq
0=\Io\be(x)e^{\f{qx}{d_I}}\big(\F{S^{\infty}\mathcal{I}_m\mathcal{I}_M}{1+me^{\f{qx}{d_I}}\mathcal{I}_m}-\F{S^{\infty}\mathcal{I}_M\mathcal{I}_m}{1+me^{\f{qx}{d_I}}\mathcal{I}_M}\big)\dx
=\Io\be(x)e^{\f{qx}{d_I}}S^{\infty}\mathcal{I}_M\mathcal{I}_m\F{me^{\f{qx}{d_I}}(\mathcal{I}_M-\mathcal{I}_m)}{(1+me^{\f{qx}{d_I}}\mathcal{I}_M)(1+me^{\f{qx}{d_I}}\mathcal{I}_m)}\dx,
\eeq
which is a contradiction.
The proof of Theorem \ref{ee ds-infty} is therefore completed. \proofend

\subsection{Concentration behaviors of EE as $d_I\rightarrow+\infty$}
In this subsection,
we present the proof of Theorem \ref{ee di-infty} as follows.

\seven

Indeed,
by a similar argument as in the proof of Theorem \ref{ee ds-infty},
we may assume that
$(S(x),I(x))\rightarrow(S_{\infty},I_{\infty})$ in $C^1([0,L])\times C^1([0,L])$,
by passing to a subsequences $d_I\rightarrow+\infty$.
Then it can be deduced that $S_{\infty},I_{\infty}\ge0$ and
\beq
\left\{\arraycolsep=1.5pt
\begin{array}{ll}
(I_{\infty})_{xx}=0, \   \ &0<x<L,\\[2mm]
(I_{\infty})_x=0,\ &x=0,L.
\end{array}
\right.
\eeq
Thus, $I_{\infty}$ is a nonnegative constant.

Assume that $I_{\infty}\equiv0$.
Applying the $S$'s equations in \eqref{equilibrium},
we then have that
$S_{\infty}(x)\equiv\s(x)$.
Set $\tilde{I}(x)=\F{I}{\|I\|_{\infty}}$.
Then we get that
\beq
\left\{\arraycolsep=1.5pt
\begin{array}{ll}
d_I \tilde{I}_{xx}-q \tilde{I}_x+(\F{\be S}{1+mI}-\gam(x))\tilde{I}=0, \   \ &0<x<L,\\[2mm]
d_I \tilde{I}_x-q \tilde{I}=0,\ &x=0,L,\\[2mm]
\|\tilde{I}\|_{\infty}=1.&
\end{array}
\right.
\eeq
By the elliptic regularity and the Sobolev embedding theorem,
by passing to a subsequence of $d_I\rightarrow+\infty$ if necessary,
we have
$\lim\limits_{d_I\rightarrow+\infty}\tilde{I}=1$ in $C^1([0,L])$.
By using \eqref{10},
there holds
\beq
\begin{split}
\Io\gam(x)\tilde{I}(x)\dx=\Io\F{\be(x)S\tilde{I}}{1+mI}\dx.
\end{split}
\eeq
This implies that when $d_I\rightarrow+\infty$,
$\Io\gam(x)\dx=\Io\be(x)\s(x)\dx$.
A contradicts to our assumption occurs. Hence $I_{\infty}>0$.
Therefore,
by \eqref{10} and the $S$'s equations in \eqref{equilibrium} again,
one can get that
\beq
\left\{\arraycolsep=1.5pt
\begin{array}{ll}
d_S (S_{\infty})_{xx}-q (S_{\infty})_x+\Lam(x)-\mu(x)S_{\infty}-\F{\be S_{\infty}I_{\infty}}{1+mI_{\infty}}+\gam(x)I_{\infty}=0,\   \ &0<x<L,\\[2mm]
d_S (S_{\infty})_x-q S_{\infty}=0,\ &x=0,L,\\[2mm]
\Io\gam(x)\dx=\Io\F{\be(x)S_{\infty}}{1+mI_{\infty}}\dx.&
\end{array}
\right.
\eeq
\proofend

\subsection{Concentration behaviors of EE as $m\rightarrow+\infty$}

\eight

Firstly, it easily follows from \eqref{equilibrium} that
$\|S(x)\|_{W^{2,p}((0,L))}\le C$,
$\|I(x)\|_{W^{2,p}((0,L))}\le C$,
since $S(x),I(x)$ are bounded for all $m\ge1$, $d_S>0,d_I>0$.
By the elliptic $L^p$ theory and Sobolev embedding theorem,
we have that
$(S(x),I(x))\rightarrow(S^*(x),I^*(x))$ in $C^1((0,L))\times C^1((0,L))$,
by passing to a subsequence of $m\rightarrow+\infty$.
Clearly, \eqref{estimate1} guarantees that $S^*(x)>0$ and $I^*\equiv0$.
Furthermore,
one can obtain that $S^*(x)=\s(x)$.

Let now $\th(x)=mI(x)$.
It is then easy to verify that
\beq
\left\{\arraycolsep=1.5pt
\begin{array}{ll}
d_I \th_{xx}-q \th_x+(\F{\be S}{1+\th}-\gam(x))\th=0, \   \ &0<x<L,\\[2mm]
d_I \th_x-q \th=0,\ &x=0,L.
\end{array}
\right.
\eeq
As the Harnack inequality yields that
$\max\limits_{x\in[0,L]}\th(x)\le C\min\limits_{x\in[0,L]}\th(x)$($C>0)$,
it indicates that
\beq
\F{1}{C}\max\limits_{x\in[0,L]}\th(x)\le\min\limits_{x\in[0,L]}\th(x)
\le\F{1}{L}\Io\th(x)\dx
=\F{1}{L}\Io mI(x)\dx\le\F{\be^*\|\Lam\|_1}{L\mu_*\gam_*}.
\eeq
Then using the standard elliptic regularity and Sobolev embedding theorem,
$\th\rightarrow\th^*\ge0$ in $C^1([0,L])$,
by passing to a subsequence $m\rightarrow+\infty$.

Suppose that $\th^*\equiv0$.
Define $\tilde{\th}=\F{\th}{\|\th\|_{\infty}}$.
Then $\|\tilde{\th}\|_{\infty}=1$ and
\beq
\left\{\arraycolsep=1.5pt
\begin{array}{ll}
d_I \tilde{\th}_{xx}-q \tilde{\th}_x+(\F{\be S}{1+\th}-\gam(x))\tilde{\th}=0, \   \ &0<x<L,\\[2mm]
d_I \tilde{\th}_x-q \tilde{\th}=0,\ &x=0,L.
\end{array}
\right.
\eeq
As before,
we may assume $\tilde{\th}\rightarrow\tilde{\th}^*\ge0$ in $C^1([0,L])$ as $m\rightarrow+\infty$.
Since $S(x)\rightarrow\s(x)$,
we then arrive at
\beq
\left\{\arraycolsep=1.5pt
\begin{array}{ll}
d_I \tilde{\th}^*_{xx}-q \tilde{\th}^*_x+(\be(x)\s(x)-\gam(x))\tilde{\th}^*=0, \   \ &0<x<L,\\[2mm]
d_I \tilde{\th}^*_x-q \tilde{\th}^*=0,\ &x=0,L.
\end{array}
\right.
\eeq
So the Harnack inequality can yield that $\tilde{\th}^*>0$.
Then $0$ is the principal eigenvalue of the eigenvalue problem \eqref{lam1(d q ga-s)},
which contradicts the assumption $\rr_0>1$.
Hence,
$\th^*>0$ on $[0,L]$ and satisfies \eqref{th*}. This proves Theorem \ref{ee m-infty}.
\proofend

\section{Discussion}
\label{discussion}

In this paper,
we investigate a reaction-diffusion-advection SIS epidemic model with saturated incidence rate and linear source.
We first establish the explicit uniform boundedness of the system \eqref{P} with $d_S=d_I$
and the implicit bounds in the general case.
Note that the
extinction and persistence of the infectious disease depend on the basic reproduction number:
the system \eqref{P} is uniformly persistent
and it admits at least one EE when the basic reproduction number greater than $1$.
Next, under a special condition \eqref{assumption},
we show in Theorem \ref{global stablity of ee} that
the EE is globally asymptotic stable.
Furthermore,
we consider the limit profiles of the basic reproduction number, as $q\rightarrow+\infty$, $d_I\rightarrow0$.
Moreover,
we study the asymptotic profiles of the EE in six cases: (i) $q\rightarrow+\infty$;
(ii) $d_S\rightarrow0$;
(iii) $d_I\rightarrow0$;
(iv) $d_S\rightarrow+\infty$;
(v) $d_I\rightarrow+\infty$;
(vi) $m\rightarrow+\infty$.

During obtaining the asymptotic behaviors of the EE,
the effects of advection, individual diffusion,
linear source and saturation are investigated.
More precisely,
Theorem \ref{ee q-infty}
shows that once the advection rate tends to infinity,
both susceptible and infected individuals will concentrate at the downstream end,
and disappear far away from the downstream end.
It follows from Theorem \ref{ee ds-0} that
as the diffusion rate of the susceptible
individuals tends to zero,
the susceptible and infected individuals should distribute inhomogeneous on the entire habitat
with a certain portion of the susceptible concentrating at the downstream end. Theorem \ref{ee di-0} indicates
that
if the diffusion rate of infected individuals approaches to zero, the susceptible distributes inhomogeneous
everywhere on the whole habitat,
and all the infected individuals will concentrate at the downstream end.
Theorem \ref{ee ds-infty} implies that upon $d_S$ tending to infinity,
the density of susceptible is positive and homogeneous on the entire habitat,
while the infectious disease always exists inhomogeneous.
Theorem \ref{ee di-infty} yields that upon $d_I$ tending to infinity,
the density of infected individuals is positive and homogeneous on the entire habitat,
while the susceptible individuals always exists inhomogeneous.
Finally,
it indicates from Theorem \ref{ee m-infty} that
when the saturated incidence rate tends to infinity,
%our result (Theorem 4.6) shows that
the EE tends to the DFE, which means that the infected individuals will vanish.

In fact,
compared with \cite{cc,sc},
there are some differences in our model:

(1).For the reaction-diffusion SIS model \eqref{P} without advection term,
when $d_S$ tends to zero,
Theorem 4.2 in \cite{sc} shows that the susceptible and infected individuals distribute inhomogeneous on the entire habitat.
However, in this case,
our result for \eqref{P} yields that besides inhomogeneous distribution,
the advection may induce the concentration
phenomenon for the susceptible individuals at the downstream end;

%(2). Our theoretical results indicate that the advection,
%especially large advection,
%may induce the concentration phenomenon.

(2) For our model \eqref{P},
the effect of the saturated incidence rate is that once $m$ tends to infinity,
the disease will be eliminated.
This indicates that the sufficiently large saturated incidence rate can help to control the disease.

\section*{Data availability}
Data sharing is not applicable to this article as obviously no datasets were generated or
analyzed during the current study.

\end{document}